\documentclass[draft,12pt,leqno]{amsart}

 \usepackage{amsmath,amssymb,amscd,amsthm}

 \usepackage{latexsym}

 \newtheorem{theorem}{Theorem}

 \newtheorem{lemma}{Lemma}
 \newtheorem{proposition}{Proposition}

 \newtheorem{corollary}{Corollary}

 \newcommand{\q}{\quad}

 \newcommand{\qq}{\quad\quad}

 \newcommand{\qqq}{\quad\quad\quad}

 \newcommand{\norm}[2]{{\left\| #1 \right\|}_{#2}}

\newcommand{\ol}[1]{\overline{#1}}

 \newcommand{\al}{\alpha}

 \newcommand{\Ga}{\Gamma}

 \newcommand{\de}{\delta}

 \newcommand{\De}{\Delta}

 \newcommand{\ve}{\varepsilon}

 \newcommand{\Th}{\Theta}

 \newcommand{\ka}{\kappa}

 \newcommand{\la}{\lambda}

\newcommand{\tf}{\tilde{f}}

\newcommand{\tg}{\tilde{g}}

 \newcommand{\La}{\Lambda}

 \newcommand{\si}{\sigma}

 \newcommand{\Om}{\Omega}

 \newcommand{\rn}{{\mathbb R}^n}

 \newcommand{\rone}{\mathbb R^1}

 \newcommand{\rtwo}{\mathbb R^2}

 \newcommand{\rthree}{\mathbb R^3}

 \newcommand{\cone}{\mathbb C^1}

 \newcommand{\stwo}{\mathbb S^2}

 \newcommand{\ca}{\mathcal A}

 \newcommand{\cc}{\mathcal C}

\newcommand{\cf}{\mathcal F}

 \newcommand{\intl}{\int\limits}

 \newcommand{\suml}{\sum\limits}

 \newcommand{\supl}{\sup\limits}

 \newcommand{\f}{\displaystyle\frac}

 \newcommand{\p}{\partial}

 \newcommand{\pp}[2]{\f{\p #1}{\p #2}}

 \newcommand{\R}{\operatorname{R}}

 \newcommand{\dpr}[2]{\langle #1, #2\rangle}

\newcommand{\ov}{\overline}
\newcommand{\lone}{{L_1^{*}}}
\newcommand{\ltwo}{{L_2^{*}}}
\newcommand{\lthree}{{L_3^{*}}}
\newcommand{\lfour}{{L_4^{*}}}
\newcommand{\hra}{\hookrightarrow}

\begin{document}

\title{On Schr\"odinger maps}

\author[Nahmod, Stefanov, Uhlenbeck]
{A. Nahmod, A. Stefanov, K. Uhlenbeck}

\address{Andrea Nahmod,  Department of Mathematics \& Statitics,
University of Massachusetts, Amherst, MA 01003, USA}

\email{nahmod@math.umass.edu}
\address{Atanas Stefanov, Department of Mathematics \&  Statitics,
University of Massachusetts, Amherst, MA 01003, USA}
\email{stefanov@math.umass.edu}
\address{Karen Uhlenbeck, Department of Mathematics, University of Texas, 
Austin, TX 78712}
\email{uhlen@math.utexas.edu}

\begin{abstract}
We study the question of well-posedness of the Cauchy problem for Schr\"odinger maps from $\rone \times \rtwo$ to the sphere $\stwo$ or to 
${\mathbb H^2}$, the hyperbolic space. The idea is to choose an 
appropriate gauge change so that the derivatives of the map will satisfy a certain nonlinear Schr\"odinger system of equations and then study 
this modified Schr\"odinger map system (MSM). 
We then prove local well posedness of the Cauchy problem for the MSM with 
minimal regularity assumptions on the data and outline a method 
to derive well posedness of the Schr\"odinger map itself from it. 
In proving well posedness of the MSM, the heart of the matter 
is resolved by considering truly quatrilinear forms of weighted 
$L^2$ functions. 
\end{abstract}

\maketitle

\date{}
\vspace{.3cm}
\section{Introduction}
\vspace{.5cm}
The harmonic map equation between two Riemannian manifolds is one of the
most studied equations in the modern geometric analysis. There are three
evolution equations which are derived from the same geometric considerations.
The best known one is the heat flow for harmonic maps, 
which was, in fact, used  By Eells and Sampson in one of the 
first papers
on harmonic maps.  This flow equation has been successfully studied 
by methods which in spirit depend on the  same geometric ideas used 
in the elliptic theory of harmonic mappings.

In  the last decade, the wave equation version, the wave map equation, 
has been studied by a number of mathematicians.  The work of 
Klainerman is probably the best known, and the recent work of  
Tao \cite{Tao1} \cite{Tao2}
is very promising.  The methods are quite different in spirit 
from the elliptic theory, and with the exception of the  classical 
work on the equation in $1+1$ dimensions and some specialized work , 
use little  in the spirit of gauge theoretic 
geometric methods.  In this paper, we 
obtain estimates (which are sufficient to give estimates  down to but 
not including the critical energy  space)  for the Schro\"edinger map 
equation in the special case from  $\rtwo$ to the sphere $\stwo$ (or to 
${\mathbb H^2}$, the hyperbolic space).   The 
historical development of the theory of this equation demonstrates 
the need for some geometric insight.

The general formulation of the Schr\"odinger equation, which we will 
not need here, arises from  writing the heat flow equation  for 
harmonic maps into a K\"ahler manifold $X$. Let 
$$
s:\rn\to X.
$$
The heat flow is then described by 
$$
\f{ds}{dt}= \nabla_s*ds.
$$
Since the K\"ahler manifold $X$ has an action by a complex structure
$J(s)$ in the tangent bundle, the Schr\"odinger map equation can be written
$$
J(s)\f{ds}{dt}=\nabla_s*ds.
$$
However the equation we are treating arises in a more natural fashion from
the Landau-Lifschitz equation for a macroscopic ferromagnetic continuum 
\cite{Papa} for $s:\rn\to\stwo$, by considering $\stwo$ as embedded into 
$\rthree$ and
$$
\f{ds}{dt}=-s \times \De s.
$$
To understand mathematically the one dimensional case, it is necessary to 
make a change of coordinates or gauge change, 
classically known as the Hasimoto 
transformation. A special gauge in the bundle $s*T(\stwo)$ is chosen in which 
the covariant derivative in the space direction is the ordinary 
differentiation $\nabla_x=d/dx$. In these coordinates, we have
$$
i\nabla_tu=\nabla_x^2u= \f{d^2u}{dx^2},
$$
where $u=ds/dx$ and $\nabla_t=d/dt+a_0$.  However, the curvature $R$ in the
image is given 
by 
$$
[\nabla_x,\nabla_t]=\f{d}{dx} a_0=R(ds/dx,ds/dt)=R(u,i \f{du}{dx}).
$$
Since the curvature of $X$ is constant at $1$, some simple K\"ahler
geometry gives
$$
R(u,i \f{du}{dx})=\f{i}{2} \f{d |u|^2}{dx}
$$
or
$$
a_0=\f{i}{2}|u|^2.
$$
In these coordinates, the equation becomes the usual integrable focusing
non-linear \\ Schr\"odinger equation. If we take ${\mathbb H^2}$ instead, 
we obtain the
defocusing case with a change of sign.

Chang, Shatah and Uhlenbeck were able to handle the one dimensional 
case for arbitrary surfaces and the radially symmetric case $n=2$ 
in the energey norm by an extension of this argument
\cite{Chang}. Our estimates follow those in spirit, although in two 
dimensions it is not possible to gauge away the derivative term completely.

We outline a proof of well-posedness in the 
coordinates we use.   The coordinates we use are not the coordinates 
of the map,  and we do not go into the technicalities of translating 
back and forth, primarily because the theory does not seem to be at 
this stage of development.

The plan of the paper is as follows. In Section \ref{sec:1} we give the 
coordinate change from a form of the Schr\"odinger map equation  to 
$\stwo$  to the form in which we are 
able to make estimates. The origin of the estimates would be totally 
mysterious without this explanation.  In Section \ref{sec:2} we state the 
fundamental estimates, which are  cubic and quintic non-linearities. The 
one nonlinearity that contains the derivative is by far the hardest one 
to handle. 
We also state the basic estimate  for inital data in 
$H^{\ve}$ (which  corresponds to $H^{1+ \ve}$ for the map). We then
include for convenience some estimates from \cite{Kenig} and \cite{Tao},
that are frequently used throughout the proof. 
The details of the proofs of estimates on various terms, which is the 
meat of the paper,  are in Sections \ref{sec:89}, \ref{sec:quintic}, 
\ref{sec:10}, \ref{sec:11}.  We use the mixed  space-time Hilbert 
spaces $X_{s,b}$ as introduced by Bourgain for the non-linear 
Schr\"odinger equation without derivative non-linearities. These are well suited
to study the low regularity behavior of these non-linear dispersive 
equations. It is interesting to note that at some stage, we need to use 
the mixed Lebesgue spaces $L_t^p L_x^q$ to handle the 
quintic nonlinearities. Our proof relies on and adapts from certain 
multilinear estimates recently obtained by Tao \cite{Tao} and 
Colliander-Delort-Kenig-Staffilani \cite{Kenig}. The authors are 
appreciative of the clarity, breadth and availability of the 
work of T. Tao, \cite{Tao}.  We note however that 
the heart of the matter is resolved by considering truly 
quatrilinear forms of weighted $L^2$ functions. 

We believe that  simliar results must hold in all dimensions, with 
$H^{\epsilon}$ replaced by $H^{n/2  -1  + \epsilon}$.   The sign of the 
curvature is not relavent to our equation, so the results hold for 
maps into the hyperbolic space ${\mathbb H^2}$ as well. 
In principle, it should not be difficult to extend the estimates to 
non-constant curvature surfaces, much as is done in the 
one-dimensional case (\cite{Chang}). 
 Since non-abelian gauge theory will be relevant for 
image manifolds of complex dimension larger than 1,  the case of 
higher dimension in the target is much more difficult.  Of course, we 
expect and hope that there are estimates which hold at the critical 
scaling regularity.  It is worth noting that for the  wave map equation,
the  spaces $X_{s,b}$, which we use  are not adequate enough to 
handle the critical case, so we are not surprised that the estimates work
down to but not at the critical case.

\vspace{.5cm}

\section{Formulation of the problem}
\label{sec:1}

\vspace{.5cm}

The Schr\"odinger map equation for $\rone\times\rn\to\stwo$ has a number 
of different descriptions, which are equivalent for smooth solutions. 
We describe this equation for all $n$ but consider in the rest of the paper 
only $n=2$. We start with a description in terms of the stereographic
 projection of $\stwo\smallsetminus\{N\} \to\cone$ where $N$ is the north pole. This is possibly the simplest for those 
unfamiliar with differential geometry. The estimates we obtain are in 
coordinate (gauge) system, which is dependent on the map, but  independent 
of any coordinate choice.

Let $s:\rn\to\cone\cup\{\infty\}=\stwo$. Then the energy of $s$ is 
$$
E(s)=\f{1}{2}\intl_{\rn}\suml_{j=1}^n
\f{\left|\p s/\p x_j\right|^2}{(1+|s|^2)^2}(dx)^n.
$$
A simple calculation shows that the Euler-Lagrange equations, or the 
equations for a harmonic map are 
\begin{eqnarray*}
& &\suml_{j=1}^n\pp{}{x_j}\left(\pp{s}{x_j}/(1+|s|^2)^2\right)+\\
& &2\suml_{j=1}^n\pp{}{x_j}\left(\left|\pp{s}{x_j}\right|^2/(1+|s|^2)^3\right)
=0.
\end{eqnarray*}
After a short computation, we find that this can be written as 
$$
\suml_j\nabla_j \pp{s}{x_j} =0.
$$
Here
$$
\nabla_j=\pp{}{x_j}-2\left(\pp{s}{x_j}\overline{s}\right)/(1+|s|^2),
$$
is the covariant derivative corresponding to the pull-back of the Levi-Civita 
connection on the tangent plane ${T}(\stwo)$ by the map $s$.

The heat equation would be 
$$
\pp{s}{t}=\suml_j {\nabla_j}\pp{s}{x_j}.
$$
The Schr\"odinger map equation is 
$$
\pp{s}{t}=\pm i\suml_j {\nabla_j}\pp{s}{x_j}.
$$
We will change gauge in this equation from the coordinate frame of the 
stereographic projection to a normalized frame, and rotate the frame to put
the pull-back covariant derivative $\nabla_j$ as near to  $\pp{}{x_j}$ 
as possible. Since we will lose track of the map $s$ during this process,
 we will need a set of consistency equations, which would be needed to 
recover the map $s$. So in addition to the equation 
\begin{equation}
\label{eq:k11}
\pp{s}{t}=i\suml_j {\nabla_j}\pp{s}{x_j}.
\end{equation}
we have two sets of consistency conditions:
\begin{equation}
\label{eq:k2}
\nabla_j\pp{s}{x_k}=\nabla_k\pp{s}{x_j}\qq \begin{array}{l}
j=0,1,\ldots,n\\
k=1,\ldots,n
\end{array}
\end{equation}
This can be computed and follows from the fact that the Levi-Civita
connection on $\stwo$ has no torsion.
\begin{equation}
\label{eq:k3}
[\nabla_j,\nabla_k]=-4 i\ Im(\overline{b_j}b_k) \qq \begin{array}{l}
j=0,1,\ldots,n\\
k=1,\ldots,n,
\end{array}
\end{equation}
where 
$$
b_j=\pp{s}{x_j}/(1+|s|^2).
$$
This is either a computation, or follows from the fact that the curvature 
of $\{\nabla_j\}$ is the pull-back of constant curvature on $\stwo$ 
by the map $s$. The appearance of $(1+|s|^2)$ is due to the fact 
that the coordinates are not (and cannot be) normalized. Note that the 
consistency conditions \eqref{eq:k2}, \eqref{eq:k3} above include the 
$t=x^0$ direction. 

We need the exsitence of a few derivatives on the map $s$ 
to classically prove the following.
\begin{theorem}
\label{theo:ak1}
Let $s:\rone\times\rn\to \stwo=\cone\cup\{\infty\}$ be a Schr\"odinger map
 of finite energy which is asymptotic to $0\in\cone$ at spatial infinity
(which can be assumed by rotation). Let 
\begin{eqnarray*}
u_j &=& (1+|s|^2)^{-1}e^{i\psi}\pp{s}{x_j}\\
D_j &=& (1+|s|^2)e^{i\psi}\circ\nabla_j\circ(1+|s|^2)e^{-i\psi} = \\
&=&\pp{}{x_j}+ia_j.
\end{eqnarray*}
Then for each $t$ there exists a unique choice of $\psi$ such that 
\begin{eqnarray}
\label{eq:k5}
\textup{div}\ a &=& 0;\qq a\sim 0 \qq \textup{at infinity}\\
\label{eq:k6}
u_0 &=& i\suml_j D_ju_j \\
\label{eq:k7}
D_j u_k &=& D_ku_j \qq\begin{array}{l}
j=0,1,\ldots,n\\
k=1,\ldots,n,
\end{array}\\
\label{eq:k8}
[D_j,D_k] &=& i\left(\pp{a_k}{x_j}-\pp{a_j}{x_k}\right)= \\
\nonumber
&=& -4i\ Im(u_j\overline{u_k})\qq\begin{array}{l}
j=0,1,\ldots,n\\
k=1,\ldots,n,
\end{array}.
\end{eqnarray}
\end{theorem}
\begin{proof}
Note that \eqref{eq:k7} and  \eqref{eq:k8} are gauge invariant equations.
 The transformation $\nabla_j\to D_j$ and $\pp{s}{x_j}\to u_j$ are the 
same gauge change. The choice of $\psi$ is possible because 
\begin{eqnarray*}
a_j &=& 2\ Im\left(s\ \overline{\pp{s}{x_j}}/(1+|s|^2)\right)-\pp{\psi}{x_j}= \\
&=& 2\ Im (\overline{b_j}s)-\pp{\psi}{x_j}.
\end{eqnarray*}
We simply choose a Hodge gauge with 
$$
\suml_{j=1}^n \pp{}{x_j}\left(2\ Im (\overline{b_j}s)-\pp{\psi}{x_j}\right)=0.
$$
If $\psi \sim 0$ at infinity, $\psi$ will be unique.
\end{proof}
{\bf Remark} Similar transformation might prove beneficial in 
the wave map problem as well.  We hope to report on that in a later paper 
\cite{Nahmod}. As we will see later, these change of variables simplifies 
{\it globally} the form of the non-linearity, which is somehow  dictated by 
geometric considerations.  Similar approach was succesfully
used by T. Tao in his work on wave maps \cite{Tao1}.

We also remark that if the map $s$ is in fact a solution in the classical Sobolev 
space $L^p_m\subset C^0$, then it is not difficult to show that 
$b\in L^p_{m-1}$. The whole point of the gauge change 
is that $\pp{a}{x_j}\in L^q_{m-1}$, where the product 
$b_j\overline{b_k}\in L^q_{m-1}$. So $\{a_j\}$ are slightly smoother 
than the $\{b_j\}$.  Of course, we will ultimately be interested in the 
mixed $L_t^pL_x^q$ and $X_{s,b}$ norms, but we will not go into details in 
this paper.

It is the set of equations in Theorem \ref{theo:ak1} which can be 
inverted to produce the map $s$. However a derived subset of these equations
form a well-posed nonlinear Schr\"odinger flow.
\begin{theorem}
\label{theo:ak2}
The following equations, which we call the ``modified Schr\"odinger map'' 
(MSM) follow from the Schr\"odinger map equation and consistency conditions
\begin{eqnarray*}
\pp{u_j}{t} &=& i\De u_j-2\suml_k a_k \pp{u_j}{x_k}-
\left(\suml_k a_k^2\right)u_j+ \\
&+& 2\ Im (\overline{u_j}u_k)u_j-i a_0u_j \qqq j=1,\ldots,n;
\end{eqnarray*}
where 
\begin{eqnarray*}
a_k &=& \suml_{l=1}^n \pp{\ka_{lk}}{x_l};\\
d\ \ka &=& 0;\\
\De \ka_{kj} &=& -4\ Im(u_k\overline{u_j}) \qq\begin{array}{l}
j=0,1,\ldots,n\\
k=1,\ldots,n,
\end{array}; \\
\De a_0 &=& -4 \suml_{j=1}^n\suml_{k=1}^n \left[\pp{}{x_k}\pp{}{x_j}
Re(u_k\overline{u_j})-\f{1}{2}\left(\pp{}{x_k}\right)^2u_j\overline{u_j}
\right].
\end{eqnarray*}
\end{theorem} 
\begin{proof}
$$ D_0\ u_j = i[D_j,D_k]u_k+i \suml_k D_k^2 u_j $$
and on the other hand 
$$ D_0\ u_j = D_ju_0=i D_j\suml_k (D_k u_k),$$ 
Here we have used $D_ku_j=D_ju_k$. In the first equation, we see the terms 
$$ i[D_j,D_k]u_k=-2\ Im(u_j\overline{u_k})u_k.$$ and since
\begin{eqnarray*}
D_k^2 &=& \f{\p^2}{\p x_k^2}+i a_k \pp{}{x_k}-a_k^2\ , \\
& &\left(\suml_k D_k^2\right) u_j = \De u_j=2i \suml_k a_k\pp{u_j}{x_k}-
\left(\suml_k a_k^2\right)u_j.
\end{eqnarray*}
Next, we use that $$
\textup{div}\ a=\sum_{j}\pp{a}{x_j}=0
$$ and that $$ \pp{a_j}{x_k}-\pp{a_k}{x_j}=-4\ Im(\overline{u_j}u_k)$$ 
to write $$ a=\textup{div}\ \ka=*d*\ \ka $$ and
$$ \De\ \ka=-4\ Im(\overline{u_j}u_k). $$ Now 
\begin{eqnarray*}
& &d\ a_0-\pp{}{t} a=4\ Im(\overline{u_k}u_0)\\
& &d*\ a=0 \\
& &\De a_0=-4\ \suml_j \pp{}{x_j}(Im\overline{u_j}u_0)\\
& &u_0=i \suml_jD_ju_j.
\end{eqnarray*} Then
\begin{eqnarray*}
Im(\overline{u_k}u_0)&=&-\suml_j Re(\ol{u_k}(D_ju_j)) \\
&=&-\suml_j \pp{}{x_j}Re(\ol{u_k}u_j)-Re(\ol{D_ju_k} u_j) \\
&=& -\suml_j \pp{}{x_j}Re(\ol{u_k}u_j)-Re(\ol{D_ku_j} u_j)\\
&=& \suml_j \pp{}{x_j}Re(\ol{u_k}u_j)-\f{1}{2}\pp{}{x_k}(|u_j|^2).
\end{eqnarray*} Hence, $$ \De \ a_0=4\suml_j\pp{}{x_j}Re(\ol{u_k}u_j)-\f{1}{2}\pp{}{x_k}(|u_j|^2).$$ \end{proof}
The modified Schr\"odinger map equation (MSM) is the $j=1,\ldots,n$ flows
for $u_j$ and the nonlinear operators defining the $a_j$'s. In this paper
we prove this equation is locally well-posed when the data is in 
$H^{\epsilon}$, $\epsilon>0$ for $n=2$.

It is not possible to go back directly from solutions of the MSM system 
to Schr\"odinger maps. In fact even in the one- dimensional case 
Chang-Shatah-Uhlenbeck (\cite{Chang})  do not attempt this. In that case we
have $a_1=0, a_0=-2u_1\ol{u_1}$. However, we sketch here a method of proving
local well posedness for the Schr\"odinger map for data in $H^{1+\epsilon}$. 

We assume that it should be possible to prove local well-posedness for data 
in $H^k$ for large $k$ for the map directly. Such solutions transform over to 
solutions of the complete (overdetermined) system. Our regularity results of 
Section \ref{sec:reg} show that the time of existence depends only on 
$\norm{u_0}{H^{\epsilon}}$ or the $H^{1+\epsilon}$ norm of the initial 
data for $s$. Moreover, we have estimates on the differences. So given 
an initial data $q$ in $H^{1+\epsilon}$, we approximate by smooth  
$q^{\al} \in H^k$, whose solutions $u^\al$ satisfy the full set of equations 
and consistency conditions. The solution produced by the well-posedness
result in Theorem \ref{theo:main} will be a limit of the solutions in 
$X_{\epsilon, 1/2+}$ and hence also satisfy the entire set of 
consistency conditions. 

As we have remarked before, the transformation formulas between 
$u\in X_{\epsilon, 1/2+}$ and the map $s$ are fairly complicated. However, we are able
to prove in this fashion that the equations \eqref{eq:k5},  \eqref{eq:k11}, 
 \eqref{eq:k2},  \eqref{eq:k3} are well-posed for initial data 
$u_0\in H^\epsilon$ via this circular route; modulo the lack of a published
proof that the Schr\"odinger map equation is locally well-posed in $H^k$ 
for large $k$.

\vspace{.5cm}

\section{The modified Schr\"odinger map system}
\label{sec:2}
 
\vspace{.5cm}

According to our reductions in the previous section, 
we consider the system of coupled nonlinear Schr\"odinger  equations
in $R^{2+1}$
\begin{equation}
\label{eq:1}
\displaystyle\left|\begin{array}{rl}
\p_t u_1 &= i\De u_1+2 
\left(\pp{\beta}{x_1}\pp{u_1}{x_2}-
\pp{\beta}{x_2}\pp{u_1}{x_1}\right)-i\al u_1-i|\nabla \beta|^2u_1\pm
Im(u_2 \overline{u_1})u_2, \\
\p_t u_2 &= i\De u_2+2\left(\pp{\beta}{x_1}\pp{u_2}{x_2}-
\pp{\beta}{x_2}\pp{u_2}{x_1}\right)-i \al u_2-|\nabla \beta|^2u_2
\pm  Im(u_1\overline{u_2})u_1\\
u_1(x,0) & = u^1_0(x),\\
u_2(x,0) &= u^2_0(x).
\end{array}
\right.
\end{equation}
where
\begin{eqnarray*}
\De \beta &=& \pm 2 Im (u_1\overline{u_2}), \\
\De \al &=& \pm \suml_{k,j=1}^2
2\left(\p_{x_k}\p_{x_j} Re(u_k\overline{u_j})-
\p_{x_k}^2|u_j|^2\right).
\end{eqnarray*}
Our main theorem asserts that the system \eqref{eq:1} is locally 
well-posed (the spaces $X_{s,b}$ are to be defined shortly).
\begin{theorem}
\label{theo:main}
For every $\ve>0$ and data $u_0\in H^{100\ve}$, 
there exists $T=T(\norm{u_0}{H^{100\ve}})$, such that the system \eqref{eq:1}
has a unique solution $u$ satisfying
$$
u\in C([0,T],H^{100\ve})\bigcap X_{100\ve, 1/2+\ve}.
$$
Moreover there exists constant $C_\ve$, independent of $u_0$ such that
\begin{equation}
\label{eq:700}
\norm{u}{X_{100\ve, 1/2+\ve}}\leq C_\ve\norm{u_0}{H^{100\ve}}.
\end{equation}
Finally, the map $u_0\backepsilon H^{100\ve}\to u\backepsilon
 C([0,T],H^{100\ve})\bigcap X_{100\ve, 1/2+\ve}$ is 
Lipschitz.
\end{theorem}
Essentially, we want to prove short time existence and
uniqueness for data in the Sobolev space $H^{100\ve}$ , 
provided that $u_1$ and $u_2$ are components of the solution and 
therefore live
in the same function spaces. From now on, we will not distinguish between
$u_1$ and $u_2$ as they come in our formulae, as we will only use their 
functional analytic properties, not the fact that they are solutions.
Occasionally, we will be  refering to the  vector $u=(u_1,u_2)$ and 
the data 
$u_0=(u^1_0, u^2_0)$. 

Next,  by the Duhamel's principle, one obtains the following equivalent
integral formulation for the system
\begin{equation}
\label{eq:3}
u(x,t)=e^{i t \De}u_0+\intl_0^t e^{i (t-\tau)\De}F_u(\tau, \cdot)d\tau,
\end{equation}
where $F$ is the nonlinearity consisting of four terms  in \eqref{eq:1}.
Introduce the (Schr\"odinger version) of the
{\it global} Bourgain spaces $X_{s,b}$ as the set of all functions $u$ with
$$
\int |\hat{u}(\xi,\tau)|^2<\tau-|\xi|^2>^{2b}<\xi>^{2s} d\xi
d\tau<\infty,
$$
where $<\xi>:=(1+|\xi|^2)^{1/2}$ and  $<\tau-|\xi|^2>:=(1+|\tau-|\xi|^2|^2)^{1/2}$. We also introduce the space $X_{s,b}^{-}$ as
$$
X_{s,b}^{-}:=\left\{u:\int |\hat{u}(\xi,\tau)|^2<\tau+|\xi|^2>^{2b}<\xi>^{2s} d\xi
d\tau<\infty\right\}.
$$
Note that the dual space to $X_{s,b}$ is $X_{-s,-b}^{-}$.
Sometimes we will refer to $s$ as the amount of elliptic smoothness in
$X_{s,b}$ and to
$b$ as the parabolic smoothness for the corresponding space.
We will also need a {\it local} version of the $X_{s,b}$ spaces, since
our solutions are local in nature. Define
$$
\norm{u}{X_{s,b}([0,T]\times\rtwo)}=\inf\left\{\norm{U}{X_{s,b}}:
U|_{[0,T]\times\rtwo)}=u\right\}.
$$
Sometimes, we will not distinguish between the local and the global spaces.
Our estimates are performed on a solution cut off in time in the global space.
Thus, they are in particular estimates on the solution in the local
space with time interval given by the lifespan of the solution.

It is well known (see for example \cite{Staffilani}) that the Schr\"odinger
semigroup $e^{i t \De}$ has certain smoothing effect on the parabolic 
derivatives. More precisely, let $\psi$ be a smooth characteristic function
of the interval $(-1,1)$ and  $1-\ve>b>1/2$ for some positive $\ve>0$. Then
\begin{equation}
\label{eq:2}
\norm{\psi(\de^{-1}t) \intl_0^t 
e^{i (t-\tau)\De} F(\tau,\cdot)d\tau}{X_{s,b}}\leq C_\ve\de^\ve 
\norm{F}{X_{s,b-1+\ve}}.
\end{equation}
The following estimate for the growth of the free solution $e^{i t \De}u_0$ 
in the Bourgain spaces is also well-known (see \cite{Staffilani}, 
Lemma 2.3.1).
\begin{lemma}
\label{le:132}
For
$s\geq 0$ and $b>1/2$, $0<\de_0\leq 1$, we have
$$
\norm{\psi(t/\de_0)e^{i t \De}u_0}{X_{s,b}}\lesssim \de_0^{(1-2b)/2}
\norm{u_0}{H^s}.
$$
\end{lemma}
The approach for solving \eqref{eq:1} is by the
method of Pickard iterations. Therefore, to prove short-time existence
(and uniqueness in certain class), one needs to show that the map
$$
\Phi(v)=\psi(t/\de_0)e^{i t \De}u_0+\psi(\de^{-1}t)\intl_0^t e^{i (t-\tau)\De}F_v(\tau, \cdot)d\tau,
$$
is a contraction on a ball in a suitable Banach space. We choose
$X_{100\ve,1/2+\ve}\times X_{100\ve,1/2+\ve} $, if our data 
$u_0\in H^{100\ve}\equiv X_{100\ve,0}$.
\vspace{.3cm}
Some remarks are in order. 
\vspace{.3cm}
\begin{itemize}
\item By scaling and dimensional analysis, it is easy to see that if 
$u(x,t)$ solves the initial value problem \eqref{eq:1}, then 
$u_\al(x,t)=\al u(\al x,\al^2t)$ solves the same with data 
$\al u_0(\al x)$. Hence the system \eqref{eq:1} is critical in 
$L^2(\rtwo)$ (i.e. the critical index is $s_c=0$). It is interesting to note
that although the term $|\nabla\beta|^2 u_j$ is essentially quintic, it scales 
and -as we shall see- behaves like a cubic nonlinearity. Indeed general quintic
nonlinearities have critical index $s_c=1/2$ in two dimensions. However our quintic is very special and it actually has $s_c=0$.
\item The terms $i\al u_k$ and
$\pm iIm(u_k\overline{u_j})u_k$  have
approximately the structure of a {\it cubic nonlinearity} $|u|^2 u$. 
These have been extensively studied  and estimates have been obtained in 
various Lebesgue and Besov spaces. Nevertheless, the methods in 
 \cite{Kenig} and \cite{Tao} allow us to  estimate in the Schr\"odinger
$X_{s,b}$ spaces. We will refer to these terms as $F_{cubic}$.
\item The terms $|\nabla \beta|^2 u_k$ are essentially {\it quintic} in $u$, 
which are in general  handled in  spaces with at least {\it half} 
derivative on the data.  However, as we shall see later, two of the
entries in the five-linear forms come with a ``missing derivatives''. This
allows us to use Sobolev embedding in conjunction with our new embedding
theorem for $X_{s,b}$ spaces to get the estimates with minimal 
smoothness assumptions.
We refer to these as $F_{quintic}$. More specifically,
$$
F_{quintic}(u_1,u_2,u_3,u_4,u_5)=\nabla \De^{-1}(u_1\overline{u_2})
\nabla \De^{-1}(u_3\overline{u_4})u_5.
$$ 
\item
For the first term, more refined analysis is needed,
since it involves derivatives of the solution. On the other hand, 
we need to  control this
expression with  virtually no regularity present, 
and to this end one needs to  exploit the ``null form'' structure. 
Observe that the ``null form'' nonlinearity is also (anti) trilinear.
We call it $F_{null}(u_1,u_2,u_3)$, where 
$$
F_{null}(u_1,u_2,u_3)=\pp{\beta}{x_1}\pp{u_3}{x_2}-
\pp{\beta}{x_2}\pp{u_3}{x_1}.
$$
\end{itemize}
Since we
do not distinguish between $u_1$ and $u_2$ in our estimates 
(we simply assume that they 
are in $X_{100\ve,1/2+\ve}$), it  will suffice to prove 
\begin{eqnarray}
\label{eq:907}
& &\norm{F_{cubic}(u_1,u_2,u_3)}{X_{100\ve,-1/2+2\ve}}\lesssim 
\norm{u_1}{X_{100\ve,1/2+\ve}}\norm{u_2}{X_{100\ve,1/2+\ve}}
\norm{u_3}{X_{100\ve,1/2+\ve}}.\\
& & \label{eq:909}
\norm{F_{null}(u_1,u_2,u_3)}{X_{100\ve,-1/2+2\ve}}
\lesssim 
\norm{u_1}{X_{100\ve,1/2+\ve}}\norm{u_2}{X_{100\ve,1/2+\ve}}
\norm{u_3}{X_{100\ve,1/2+\ve}}.
\end{eqnarray}
and 
\begin{equation}
\label{eq:911}
\norm{F_{quintic}(u_1,u_2,u_3,u_4,u_5)}{X_{100\ve,1/2+\ve}}\lesssim 
\prod_{j=1}^5\norm{u_j}{X_{100\ve,1/2+\ve}}
\end{equation}
Then  \eqref{eq:2}, \eqref{eq:907}, \eqref{eq:909},\eqref{eq:911}  
and Lemma \ref{le:132} imply that  
$$
\Phi:B_{R_0}(X_{100\ve,1/2+\ve})\times B_{R_0}
(X_{100\ve,1/2+\ve})\to B_{R_0}(X_{100\ve,1/2+\ve})\times 
B_{R_0}(X_{100\ve,1/2+\ve})
$$
is a contraction mapping for suitably chosen $R_0$ and $\de$.

Indeed, our estimates for the nonlinearities 
 will suffice to show that for $v=(v_1,v_2)$
\begin{eqnarray*}
\norm{\Phi(v)}{X_{100\ve,1/2+\ve}\times X_{100\ve,1/2+\ve}}&\lesssim &
\norm{u_0}{H^{100\ve}}+\de^\ve\prod\limits_{j=1}^3\norm{v_j}
{X_{100\ve,1/2+\ve}}+\\
&+& \de^\ve\prod\limits_{j=1}^5\norm{v_j}{X_{100\ve,1/2+\ve}},
\end{eqnarray*}
where $\de_0=1$ from Lemma 
\ref{le:132} and $v_j=v_1$ or $v_2$ for $j=3,4,5$. 
Hence if $R_0\sim \norm{u_0}{H^{100\ve}}$ and $0<\de\leq \de_0$
satisfies 
$$
\de<\min\left(\f{1}{\norm{u_0}{H^{100\ve}}^{2/\ve}},
\f{1}{\norm{u_0}{H^{100\ve}}^{4/\ve}}\right),
$$
we have that 
$$
\Phi:B_{R_0}\times B_{R_0}
\to B_{R_0}\times 
B_{R_0}
$$
and that there exists $0< C_0 <1$ such that 
$$
\norm{\Phi(v)-\Phi(\tilde{v})}
{X_{100\ve,1/2+\ve}\times X_{100\ve,1/2+\ve}}\leq C_0\norm{v-\tilde{v}}
{X_{100\ve,1/2+\ve}\times X_{100\ve,1/2+\ve}}.
$$ In particular we have estimates on the differences. 
The above is also enough to establish the Lipschitz bounds that were claimed in 
Theorem  \ref{theo:main}.

The following 
Lemma yields the estimates needed to handle the cubic-like 
nonlinearities.
\begin{lemma}
\label{le:956}
\begin{eqnarray}
\label{eq:371}
& &\norm{u_1\overline{u_2}u_3}{X_{s,-1/2+2\ve}}
\lesssim \norm{u_1}{X_{s,1/2+\ve}}
\norm{u_2}{X_{s,1/2+\ve}}\norm{u_3}{X_{s,1/2+\ve}} \\
& &\norm{u_1\overline{u_2} \overline{u_3}}{X_{s,-1/2+2\ve}}
\lesssim \norm{u_1}{X_{s,1/2+\ve}}
\norm{u_2}{X_{s,1/2+\ve}}\norm{u_3}{X_{s,1/2+\ve}}
\label{eq:372} \\
& & \norm{u_1 u_2 u_3}{X_{s,-1/2+2\ve}}
\lesssim \norm{u_1}{X_{s,1/2+\ve}}
\norm{u_2}{X_{s,1/2+\ve}}\norm{u_3}{X_{s,1/2+\ve}}
\label{eq:377}
\end{eqnarray}
provided $s>5\ve$.
\end{lemma}
We will prove Lemma \ref{le:956}, together with \eqref{eq:907}
for the cubic nonlinearities in Section \ref{sec:89}. The mixed space-time 
 Lebesgue
spaces are defined as the set of all functions $u$ with
$$
\norm{u}{L_t^pL_x^q}=\left(\int\left(\int |u(x,t)|^q dx\right)^{p/q}dt\right)^{1/p}.
$$
The next lemma supplies an embedding theorem for $X_{s,b}$ spaces 
into a mixed norm  Lebesgue spaces\footnote{The embedding (\ref{eq:876})
   actually holds in more generality into the $L^q_t L^r_x$ spaces
  with $2/q + n/r =n/2$ and $n$ the space dimension. The proof of it 
relies on the Strichartz inequalities and it was pointed out  
to us by T. Tao after we had derived the embedding -in our
  restricted range- as a consequence of our inequalities (\ref{eq:6754}) 
and (\ref{eq:6755}). Since we actually need  (\ref{eq:6754}) 
and (\ref{eq:6755}) per se to treat the quintic nonlinearity we prefer to keep 
the statement of such result in the fashion of Lemma \ref{le:956} above.} 
\begin{lemma}
\label{le:456}
For $1\leq p\leq 2$, 
\begin{equation}
\label{eq:876}
X_{0,1/2+}\hra L_t^{2p'}L_x^{2p}(\rtwo\times \rone)
\end{equation}

More generally, we have the bilinear estimates
\begin{eqnarray}
\label{eq:6754}
\norm{u v}{L_t^{p'}L_x^{p}}\lesssim 
\norm{u}{X_{0,1/2+\ve}}\norm{v}{X_{0,1/2+\ve}}.\\
\label{eq:6755}
\norm{u \overline{v}}{L_t^{p'}L_x^{p}}\lesssim 
\norm{u}{X_{0+,1/2+\ve}}\norm{v}{X_{0+,1/2+\ve}}
\end{eqnarray}
\end{lemma}
We prove Lemma \ref{le:456} as well as the estimates required
for the quintic nonlinearities in Section \ref{sec:quintic}.

Next, note that the ``null form'' nonlinearity $F_{null}$ has two components
 (one for each equation in the system), but both
components  look identical except for the dependence upon $u_1$ and $u_2$,
which is irrelevant in our argument. 
Thus, it suffices to consider only 
one component. 
We test $F_{null}$ 
against a function $W\in X_{-100\ve, 1/2-2\ve}^{-}$ to get
the (anti) multilinear form
$$
M(u_3,u_1,u_2,W)=
\int \left(\pp{\beta}{x_1}\pp{u_3}{x_2}-\pp{\beta}{x_2}\pp{u_3}{x_1}\right)
\overline{W} dxdt.
$$
Hence, the following theorem takes care of the null-form nonlinearity.
\begin{theorem}
\label{theo:10}
\begin{equation}
\label{eq:5}
|M(u_1,u_2,u_3,W)|\lesssim \norm{u_3}{X_{100\ve,1/2+\ve}}
\norm{u_1}{X_{100\ve,1/2+\ve}}\norm {u_2}{X_{100\ve,1/2+\ve}}
\norm{\overline{W}}{X_{-100\ve, 1/2-2\ve}},
\end{equation}
where $\De \beta= \pm 2Im(u_1\overline{u_2})$. 
\end{theorem}
We make some reductions for the proof of Theorem \ref{theo:10}. 
First, we will assume that all four functions $u_3, u_1, u_2, W$
are test functions with norm one in the corresponding spaces. Standard 
approximation techniques will then yield  
the general result.

Next, we perform integration by parts in the definition of $M$ to take
the derivatives off  $\beta$. The special cancelation properties of the expression and the lack of boundary terms allow us to rewrite $M$ as
$$
M(u_3,u_1,u_2,W)=\int  \beta
\left(\pp{\ov{W}}{x_2}\pp{u_3}{x_1}-\pp{\overline{W}}{x_1}\pp{u_3}{x_2}\right).
$$
Since $\De \beta=\pm 2Im(u_1\overline{u_2})=
\pm (u_1\overline{u_2}-u_2\overline{u_1})$, two similar terms arise. 
By slightly abusing our notations, we will call one of them 
(say the one that corresponds to $u_1\overline{u_2}$) $M$. 
We have
\begin{equation}
\label{eq:15}
M(u_3,u_1,u_2,W)=\int u_1 \overline{u_2} G,
\end{equation}
where
$$G=\De^{-1}\left(\pp{\overline{W}}{x_2}\pp{u_3}{x_1}-\pp{\ov{W}}{x_1}
\pp{u_3}{x_2}\right).$$
Parseval's identity, together with
$\displaystyle\overline{\widehat{u_2}}(\xi,\tau)=
\widehat{\overline{u_2}}(-\xi,-\tau)$ imply that
\begin{equation}
M(u_3,u_1,u_2,W)=\int u_1G \overline{u_2}=\int \widehat{u_1 G}(\xi,\tau)
\widehat{\overline{u_2}}(-\xi,-\tau)d\xi d\tau.
\end{equation}
Since the complex conjugation is an isometry of $X_{s,b}$ onto $X_{s,b}^{-}$,
we can write $$
\widehat{\overline{u_2}}(-\xi,-\tau)=
\f{h_2(-\xi,-\tau)}{<\tau-|\xi|^2>^{1/2+\ve}}<\xi>^{-100\ve},$$ for
some $h_2\in L^2$.

Similarly by  using  the properties of $X_{s,b}$ spaces, we express
$u_1(u_2)$ via its $L^2$ representative $h_1$ ($h_2$ respectively)
 and some weights dictated by the
particular space. We get,
\begin{eqnarray}
\label{eq:8}
M(u_3,u_1,u_2,W) &=& \int
\f{h_1(\xi-\eta,\tau-\mu)<\xi-\eta>^{-100\ve}}{<\tau-\mu-|\xi-\eta|^2>^{1/2+\ve}}
\widehat{G}(\eta,\mu) \times\\
&\times &\f{h_2(-\xi,-\tau)<\xi>^{-100\ve}}{<\tau-|\xi|^2>^{1/2+\ve}}
d\xi d\eta d\tau d\mu,
\nonumber
\end{eqnarray}
where
\begin{eqnarray}
\label{eq:9}
 \widehat{G}(\eta,\mu)  &=& \f{1}{|\eta|^2} \int (\eta_1 z_2-\eta_2 z_1)
\widehat{\ov{W}}(\eta-z,\mu-s)
\widehat{u_3}(z,s)dz ds=\\
&=& \f{1}{|\eta|^2} \int \dpr{\eta}{z^{\bot}}
\widehat{\ov{W}}(\eta-z,\mu-s)
\widehat{u_3}(z,s)dz ds.
\nonumber
\end{eqnarray}
We will need the representation
of $M$ as a quatrilinear form applied to four $L^2$ functions. Meanwhile, we
change variables $\eta\to\eta+z$ and $z\to -z$ to  obtain
\begin{eqnarray}
\label{eq:100}
\La(h_1,h_2,f,g) &=& \int
\f{h_1(\xi-\eta+z,\tau-\mu)<\xi-\eta+z>^{-100\ve}}{<\tau-\mu-|\xi-\eta+z|^2>^
{1/2+\ve}} \f{h_2(-\xi,-\tau)<\xi>^{-100\ve}}{<\tau-|\xi|^2>^{1/2+\ve}} \\
&\times & \f{\dpr{\eta}{z^{\bot}}}{|\eta-z|^2}
\f{f(\eta,\mu-s)<\eta>^{100\ve}}{<\mu-s-|\eta|^2>^{1/2-2\ve}}
\f{g(-z,s)<z>^{-100\ve}}{<s-|z|^2>^{1/2+\ve}}
d\xi d\eta dz d\tau d\mu ds.
\nonumber
\end{eqnarray}

 We break up the $\eta$ and $z$ integration in the definition of
$\La$ to obtain
$$
\La(h_1,h_2,f,g)=\intl_{|z|\leq |\eta|/2\q \&\q |\eta|\leq|z|/2 } \ldots  +
\intl_{|z|\sim |\eta|} \ldots =
\La_{off diag}+\La_{diag}
$$
We will estimate  $\La_{off diag}$ in Section \ref{sec:10}  and
$\La_{diag}$ in Section \ref{sec:11}.

\vspace{.3cm}

\section{Some remarks regarding multilinear forms}

\vspace{.3cm}

In this section, we  follow \cite{Tao} to introduce somewhat
 more general framework for
the multilinear forms that we have to deal with.

For an integer $k$ and abelian group $Z$ , define
the hyperplane $$\Ga_k(Z)= \{(\xi_1,\ldots, \xi_k)\in Z^k:
\xi_1+\ldots+\xi_k=0\}.$$
A $[k,Z]$ multiplier is a function $m:\Ga_k(Z)\to \cc$, so that there exists a constant $C$, such
that the inequality
$$
|\intl_{\Ga_k(Z)}m(\xi)\prod\limits_{j=1}^k f_j(\xi_j)|\leq C \prod\limits_{j=1}^k \norm{f_j}{L^2(Z)},
$$
holds for all $f_1, \ldots f_k\in L^2(Z)$. The best constant $C$ with the above property is
naturally called a multiplier norm for $m$ and is usually denoted $\norm{m}{[k,Z]}$. We will also need
the notation $\Ga_k(Z,\xi_j=\eta_j):=\Ga_k(Z)\cap \{\xi_j=\eta_j\}$ for a fixed  $\eta_j$.

The Cauchy-Schwarz inequality gives the following very useful lemma. (cf. \cite{Tao}, Lemma 3.9)
\begin{lemma}
\label{le:1}
If $m$ is a $[k,Z]$ multiplier, then
$$
\norm{m}{[k,Z]}\leq \supl_{\eta_j\in Z}\left(\intl_{\Ga_k(Z,\xi_j=\eta_j)}|m(\xi)|^2\right)^{1/2}.
$$
\end{lemma}
We will also need the following corollary for the cases $k=3,4$ (cf. \cite{Tao}, Corollary 3.10).
\begin{corollary}
\label{cor:1}
For any subsets $A, B, C$ of $Z$, we have
\begin{eqnarray}
\label{eq:10}
& &\norm{\chi_A(\xi_1)\chi_B(\xi_2)}{[3,Z]}\leq |\{\xi_1\in A: \xi-\xi_1\in B\}|^{1/2}, \\
& & \norm{\chi_A(\xi_1)\chi_B(\xi_2)\chi_C(\xi_3)}{[4,Z]}\leq
|\{(\xi_1,\xi_2)\in A\times B: \xi-\xi_1-\xi_2\in C\}|^{1/2}
\label{eq:11}
\end{eqnarray}
for some $\xi\in Z$.
\end{corollary}
\begin{proof}
The proofs of \eqref{eq:10} and  \eqref{eq:11} are similar, so 
we show only \eqref{eq:11}.
By Lemma \ref{le:1},
\begin{eqnarray*}
& &\norm{\chi_A(\xi_1)\chi_B(\xi_2)\chi_C(\xi_3)}{[4,Z]}\leq
\supl_{\eta_4\in Z}\left(\intl_{\Ga_k(Z,\xi_4=\eta_4)}|\chi_A(\xi_1)|^2
|\chi_B(\xi_2)|^2|\chi_C(\xi_3)|^2\right)^{1/2} \lesssim  \\
&\lesssim&\norm{|\chi_A\otimes \chi_B|*|\chi_C|}{L^\infty(Z)}^{1/2}=
|\{(\xi_1,\xi_2)\in A\times B: \xi-\xi_1-\xi_2\in C\}|^{1/2}.
\end{eqnarray*}
\end{proof}
The other important technical tool that is crucial for us
will be a form of the Schur's test and the  box
localization method (cf. Lemma 3.11 and Corollary 3.13 in \cite {Tao}).
We define the $j$'th support
of $m$ to be the set
$$
\textup{supp}_j(m)=\{\eta_j\in Z:\Ga_k(Z,\xi_j=\eta_j)\cap \textup{supp}(m)\neq\emptyset\}.
$$
In particular, if $\displaystyle m(\xi_1,\ldots, \xi_s)=\prod\limits_{j=1}^s m_j(\xi_j)$, we have
$\textup{supp}_j(m)=\textup{supp}(m_j)$.
\begin{lemma}
\label{le:2}
Let $R$ be a rectangular box in $Z$. Suppose also that $\textup{supp}_j(m)$ is contained in a $R+\eta_j$ for some
$\eta_j\in Z $ and $1\leq j\leq k-2$. Then
$$
\norm{m}{[k,Z]}\sim \supl_{\eta_{k-1},\eta_k\in Z} \norm{m\chi_{R+\eta_{k-1}}(\xi_{k-1})
\chi_{R+\eta_{k}}(\xi_k)}{[k,Z]}.
$$
\end{lemma}
In particular, the box localization principle states, that if we have a multiplier in which all
but two of the variables are restricted to  sets of certain diameter, we can restrict (at the expense
of a bigger constant) the remaining
two variables to sets with the same diameter.

Next, consider a fixed smooth function  $\psi$ on $\rone$ supported
around $1$. Introduce 
also a cutoff $\psi_0$  on $\rone$ supported around $0$ with 
$$
\psi_0(\cdot)+\suml_{R\geq 1} \psi(R^{-1}(\cdot))\equiv 1,
\qq \textup{ on } \q \rone.
$$
Let $\widehat{u_R}=\widehat{S_R u}=
\psi(R^{-1}|\cdot|)\widehat{u}(\cdot,t)$ be the standard Littlewood-Paley operator {\it in space}  
at frequency
$R$ applied to the function $u(\cdot,t)$. Sometimes we will slightly
abuse notations  by using $f_R$ to denote the restriction of the function
$f(\xi)$ to the annulus $\{\xi: |\xi|\sim R\}$.
We also   introduce  the
following notation: for every sequence of real numbers $L_1,\ldots, L_n$,
the sequence $L_1^*, \ldots, L_n^* $  will denote the permutation in
 {\it increasing} order of the original sequence.

We need the following lemmas, which are adaptation of estimates  
$(80)$ and $(86)-(89)$ in \cite{Tao} 
(cf. Proposition $11.1$ and Proposition 11.2 in \cite{Tao}).
\begin{lemma}
\label{le:5}
For the $[3,R^{d+1}]$ multiplier
$m_1(\xi,\tau)=
\prod_{i=1}^3\chi_{ \tau_i-|\xi_i|^2\sim L_i}
\chi_{|\xi_i|\sim R_i}$, we have the estimate 
$$
\norm{m_1}{[3,R^{d+1}]}\lesssim \lone^{1/2}{R_3^*}^{-1/2}
{R_1^*}^{(d-1)/2}\min(R_1^* R_3^*,\ltwo)^{1/2}. 
$$
\end{lemma}
\begin{lemma}
\label{le:11}
For $\xi_0\in {\mathbb R^d}$ and the multiplier 
$$
m_2(\xi,\tau)=
\chi_{ \tau_1-|\xi_1|^2\sim L_1}\chi_{ \tau_2+|\xi_2|^2\sim L_2}
\chi_{|\xi_1-\xi_0|\leq r,|\xi_2+\xi_0|\leq r,|\xi_3|\leq r }
$$
we have
\begin{itemize}
\item If $|\xi_0|\lesssim r$, then 
\begin{equation}
\label{eq:523}
\norm{m_2}{[3,R^{d+1}]}\lesssim \min(L_1,L_2)^{1/2}r^{(d-2)/2}\min(r^2, \max(L_1,L_2))^{1/2}.
\end{equation}
\item If $|\xi_0|\gg r$ and $H=||\xi_1|^2+|\xi_3|^2-|\xi_2|^2|\sim \tau_3-|\xi_3|^2, H\lesssim |\xi_0|r$ then
\begin{equation}
\label{eq:524}
\norm{m_2}{[3,R^{d+1}]}\lesssim \min(L_1,L_2)^{1/2}\f{r^{(d-1)/2}}{|\xi_0|^{1/2}}\min\left(H, \f{H}{r^2} \max(L_1,L_2)\right)^{1/2}.
\end{equation}

\item If $|\xi_0|\gg r$ and $H\nsim \tau_3-|\xi_3|^2, H\lesssim |\xi_0|r$  then
\begin{equation}
\label{eq:525}
\norm{m_2}{[3,R^{d+1}]}\lesssim \min(L_1,L_2)^{1/2}\f{r^{(d-1)/2}}{|\xi_0|^{1/2}}\min(|\xi_0|r, \max(L_1,L_2))^{1/2}.
\end{equation}
\end{itemize}
\end{lemma}
\begin{proof}
The proof is a reprise of the argument behind Proposition $11.2$ in 
\cite{Tao}, so we will just indicate the main points.
\begin{itemize}
\item If $|\xi_0|\lesssim r$, we imagine that we
 have the additional restriction $\tau_3-|\xi_3|^2\sim L$ in 
the multiplier $m_2$, which will be artificial and it 
will not play any role. Compute $|H|=||\xi_1|^2+
|\xi_3|^2-|\xi_2|^2|\lesssim r^2$. Thus, we might be in any
of the cases $(86)$, $(88)$, $(89)$ in \cite{Tao}, but 
in all of them, we get
$$
\norm{m_2}{[3,R^{d+1}]}\lesssim \min(L_1,L_2)^{1/2}r^{(d-2)/2}
\min(r^2, \max(L_1,L_2))^{1/2}.
$$
\item If $|\xi_0|\gg r$, we will impose again the additional artificial restriction $\tau_3-|\xi_3|^2\sim L$. Then 
$|H|=||\xi_1|^2+|\xi_3|^2-|\xi_2|^2|\lesssim |\xi_0|r$ and the rest follows from Proposition $11.2$ in \cite{Tao}.
\end{itemize}
\end{proof}

{\bf Remark} Geometric considerations indicate that the restriction 
$\tau_k\pm|\xi_k|^2\sim L_k \sim \lthree$ 
is very weak or redundant altogether. Actually, in the proof of
Lemma  \ref{le:5} and Lemma \ref{le:11} (see the discussion 
and the reductions in \cite{Tao}, estimates $(33)-(40)$), one always
estimates by the norm of the  multiplier, where the restriction 
$\tau_k\pm|\xi_k|^2\sim L_k 
\sim \lthree$ is not present. Later on, when we need to break up the 
integrals in the multilinear forms  relative to the size of the weights 
$\tau_i\pm|\xi_i|^2$, we will implicitly use the fact that the restriction 
$\tau_k\pm|\xi_k|^2\sim L_k$ ($L_k\sim\lthree$) does not appear.

The  lemma below  appears in \cite{Kenig} and essentially follows by
combining various cases in  Proposition $11.2$ of \cite{Tao}. We state it
separately, since it will be used in this form in the sequel.
\begin{lemma}(Estimate 2.20 in \cite{Kenig})
\label{le:3}
For the trilinear form
$$
C(f,g,h)=\intl_{\Ga_3(R^{2+1})} f_R(\xi_1,\tau_1) \f{g_M(\xi_2,\tau_2)}{<\tau_1-|\xi_1|^2>^{1/2+\ve}}
\f{h_N(\xi_3,\tau_3)}{<\tau_2+|\xi_2|^2>^{1/2+\ve}}
$$
there is the estimate
\begin{equation}
\label{eq:6432}
|C(f,g,h)|\lesssim \left(\f{\min(M,N))}{\max(M,N)}\right)^{1/2}\norm{f_R}{L^2(\rthree)}
\norm{g_M}{L^2(\rthree)}
\norm{h_N}{L^2(\rthree)}.
\end{equation}
\end{lemma}
As a corollary, one has the estimate for products
\begin{eqnarray}
\label{eq:56}
& &\norm{S_R( (u_1)_M(\overline{u_2})_N)}{L^2(\rthree)}\lesssim
\left(\f{\min(M,N))}{\max(M,N)}\right)^{1/2} \max(M,N)^{-s}
\norm{u_1}{X_{s,1/2+\ve}}
\norm{u_2}{X_{s,1/2+\ve}} \\
& &\lesssim \left(\f{\min(M,N))}
{\max(M,N)}\right)^{1/2} R^{-s}\norm{u_1}{X_{s,1/2+\ve}}
\norm{u_2}{X_{s,1/2+\ve}}.
\nonumber
\end{eqnarray}
We have used $\max(M,N)\gtrsim R$, since otherwise $\textup{supp}
\q \cf((u_1)_M (\ov{u_2})_N) \subset \{\xi:|\xi|\ll R\}$ and
$S_R((u_1)_M (\ov{u_2})_N)=0$.

\vspace{.5cm}

\section{Cubic nonlinearities}
\label{sec:89}

\vspace{.5cm}

We start off with a proposition, that is essentially equivalent to 
Lemma \ref{le:956}.
Later on, we  will also use it in the {\it off-diagonal} 
considerations for the 
``null form'' nonlinearity.
\begin{proposition}
\label{prop:lpk}
Let $m$ be a bounded function. Then for every $\ka:1/2>\ka>5\ve$, 
the quatrilinear form 
\begin{eqnarray*}
H(f,g,h,w)&=&\intl_{\begin{array}{l}\xi_1+\xi_2+\xi_3+\xi_4=0, \\
\tau_1+\tau_2+\tau_3+\tau_4=0
\end{array}}
m(\xi_i) \f{f(\xi_1,\tau_1)<\xi_1>^{-\ka}}{<\tau_1-|\xi_1|^2>^
{1/2+\ve}} \f{g(\xi_2,\tau_2)<\xi_2>^{-\ka}}{<\tau_2\pm |\xi_2|^2>^{1/2+\ve}} \\
&\times & 
\f{h(\xi_3,\tau_3)<\xi_3>^{\ka}}{<\tau_3-|\xi_3|>^{1/2-2\ve}}
\f{w(\xi_4,\tau_4)<\xi_4>^{-\ka}}{<\tau_4\pm|\xi_4|>^{1/2+\ve}}
d\xi_1\ldots d\xi_4 d\tau_1\ldots d\tau_4.
\end{eqnarray*}
satisfies
$$
|H(f,g,h,w)|\lesssim \norm{m}{\infty}
\norm{f}{2}\norm{g}{2}\norm{h}{2}\norm{w}{2}.
$$
\end{proposition}

{\bf Remark} We will need to have a nontrivial $m$ to cover some cases
where we have a {\it zero} order pseudodifferential operators acting on
some of the entries.
\begin{proof}
Observe that one can assume $m=1$ without loss of generality
since all the proofs proceed by putting  absolute values inside the 
integrals anyway.

For the case of weight $\tau_4-|\xi_4|^2$ and $\tau_2+|\xi_2|^2$ write 
\begin{eqnarray*}
& &\tilde{f}(\xi,\tau)=
\f{f(\xi,\tau)<\xi>^{-\ka}}{<\tau-|\xi|^2>^{1/2+\ve}}\\
& &\tilde{g}(\xi,\tau)=
\f{g(\xi,\tau)<\xi>^{-\ka}}
{<\tau+|\xi|^2>^{1/2+\ve}}
\end{eqnarray*}
By taking a dyadic decomposition on the space frequency, we have 
$$
H=\suml_R\int S_R (\tilde{f}\tilde{g}) S_R(\tilde{G})=\suml_{M,N,R}
\int S_R(\tf_M \tg_N) \tilde{G}_R d\eta d\mu,
$$
where
$$
\widehat{\tilde{G}}(\eta,\mu)=\int\f{h(\eta-z,\mu-s)
<\eta-z>^{\ka}}
{<\mu-s-|\eta-z|^2>^{1/2-2\ve}}
 \f{w(z,s)}{<z>^{\ka}
<s-|z|^2>^{1/2+\ve}}dz ds.
$$
An application
of Cauchy-Schwarz and \eqref{eq:56} from Lemma \ref{le:3}  yield
$$
|H|\leq \suml_{M,N,R}\norm{S_R(\tf_M \tg_N)}{L^2}
 \norm{\tilde{G}_R}{L^2}\lesssim 
\suml_R R^{-\ka} \norm{f}{2}\norm{g}{2}\norm{\tilde{G}_R}{L^2}.
$$

To complete the desired estimate  it will suffice to show that
\begin{equation}
\label{eq:17}
\sum_{R>1} R^{-\ka}\norm{\tilde{G}_R}{L^2}\lesssim 1.
\end{equation}

We test $\widehat{\tilde{G}_R}$ against an unimodular $L^2$ function $v$
 to get the trilinear form $S$

\begin{eqnarray}
\label{eq:2317}
S(h,w,v)&:=& \dpr{\widehat{\tilde{G}_R}}{v}=\suml_{r, R} \int 
\f{h_R(\eta-z,\mu-s)<\eta-z>^{\ka}}
{<\mu-s-|\eta-z|^2>^{1/2-2\ve}}\times \\
\nonumber
&\times& \f{w_r(z,s)}{<z>^{\ka}
<s-|z|^2>^{1/2+\ve}} v_{\max(r,R)}(\eta,\mu)d\eta dz
d\mu ds
\end{eqnarray}
applied to the $L^2$ functions $h,w,v$. We estimate away the elliptic
 weights $<\eta-z>^{\ka}<z>^{-\ka}$ to get 
\begin{eqnarray*}
S(h,w,v)&\lesssim &=\suml_{r\leq R}
\left(\f{R}{r}\right)^{\ka} \int 
\f{h_R(\eta-z,\mu-s)}
{<\mu-s-|\eta-z|^2>^{1/2-2\ve}}\times \\
&\times& \f{w_r(z,s)}
{<s-|z|^2>^{1/2+\ve}} v_R(\eta,\mu)d\eta dz
d\mu ds.
\end{eqnarray*}
for positive $h, w, v$. 
Let  $L_1=\mu-s-|\eta-z|^2$, $L_2=s-|z|^2$ and $L_3=\mu-|\eta|^2$.
An application  of Lemma \ref{le:5}  yields
\begin{eqnarray*}
& &|S(h,w,v)| \lesssim \suml_{r\leq R}\suml_{\lone,\ltwo}
\f{\lone^{1/2}\min(\ltwo,Rr)^{1/2}}{<L_1>^{1/2-2\ve}<L_2>^{1/2+\ve}}
\norm{h_R}{}\norm{w_r}{}\norm{v_R}{}\left(\f{r}{R}\right)^{1/2-\ka}\lesssim\\
&\lesssim& \suml_{r\leq R}\suml_{\ltwo<Rr}\ltwo^{2\ve}
\norm{h_R}{}\norm{w_r}{}\norm{v_R}{}
\left(\f{r}{R}\right)^{1/2-\ka}+ \\
&+& \suml_{r\leq R}\suml_{\ltwo>Rr} \f{(Rr)^{1/2}}{\ltwo^{1/2-2\ve}}
\norm{h_R}{}\norm{w_r}{}\norm{v_R}{}
\left(\f{r}{R}\right)^{1/2-\ka}\lesssim
R^{4\ve} \norm{h}{}\norm{w}{}\norm{v}{}.
\end{eqnarray*}
Thus $\norm{\tilde{G}_R}{L^2}\lesssim R^{4\ve}$ and therefore \eqref{eq:17} holds.

For the  case of weight $\tau_4+|\xi|^2$, $\tau_2+|\xi_2|^2$ we need to consider 
different  pairing of our functions: $f ,h$ versus $g, w$ rather 
than $f, g$ versus $h, w$ as we have just done. One obtains two trilinear 
forms where the weights have {\it the same} signs and we apply Lemma 
\ref{le:5} to each one of them. 
Then we  perform similar and in fact simpler argument to the one 
presented above. 
In the case $\tau_2-|\xi_2|^2$,  $\tau_4-|\xi_4|^2$
we once again pair functions with {\it same signs} 
weights   and use Lemma \ref{le:5}. 
Finally, in the case $\tau_4+|\xi_4|^2$, $\tau_2-|\xi_2|^2$  the argument is 
identical to the one presented above. We omit the details.
\end{proof}

As we have already mentioned, Proposition \ref{prop:lpk} implies 
Lemma \ref{le:956}. We will show now \eqref{eq:907} for all cubic-like
nonlinearities.
We have from the defining relations for $\al$
$$
\al_j=i \suml_{j,k=1}^2R_jR_k \ Re(u_k \overline{u_j}) -
i|u_1|^2-i|u_2|^2,
$$
where $R_j$ is the Riesz transform in the $j$th variable. Since 
the multiplier corresponding to 
$R_j$ is $\xi_j/|\xi|$, \eqref{eq:907} for the nonlinearity
$\al u_j$ reduces to Proposition \ref{prop:lpk} 
with a suitable choice of $m$. 

\vspace{.5cm}
\section{Quintic nonlinearities}
\label{sec:quintic}

\vspace{.5cm}

In this section, we will estimate the terms with a ``quintic'' nonlinearities.
As it was mentioned earlier, quintic nonlinearities are difficult to
control in a space with less than a half derivative. 
We have however a very special form of the nonlinearity, which makes 
it tractable. A good model expression of what we are dealing with is 
$$
P_{-1}(u_1\overline{u_2})P_{-1}(u_3 \overline{u_4})u_5,
$$
where $P_{-1}$ will be a smoothing pseudodifferential operator of order $-1$.
At the first step, we use the cubic estimates outlined in Section \ref{sec:89}
and then we exploit the ``smoothing'' provided by $P_{-1}$ via 
Sobolev embedding.
To accomplish this program,  we will need to pass from $X_{s,b}$ to a mixed
norm Lebesgue spaces and back. It is also a question of independent interest
to study the relation between the $X_{s,b}$ spaces and the 
mixed norm Lebesgue spaces. Actually, the reader 
is probably aware that most of the current existence 
results are in fact proved via 
an appropriate contraction mapping argument in mixed Lebesgue spaces of
various sorts. We also point out that the case $p=2$ in Lemma \ref{le:456}
is the well-known Bourgain's lemma (\cite{Bourgain}, Corollary $3.39$)
\begin{proof}(Lemma \ref{le:456}).
First, we note that the bilinear estimate \eqref{eq:6754} 
implies the embedding 
\eqref{eq:876}, if one takes $u=v$. 
As it has been already noted, the endpoint $p=2$ is contained in 
Corollary $3.39$ of \cite{Bourgain},
but can also be obtained from Lemma \ref{le:5}. 
For the bilinear estimate \eqref{eq:6755}, one uses as an endpoint 
$L^2$ result \eqref{eq:6432}. Note that one needs a little bit of extra
regularity in $u, v$ ($X_{0+,1/2+}$) in order to be able to sum 
\eqref{eq:6755} in $R$. 

By complex interpolation, it 
will suffice to show  the other endpoint $p=1$. The proof for both
\eqref{eq:6754} and \eqref{eq:6755} is the same for $p=1$, so we concentrate
on \eqref{eq:6754}. 
Since $\norm{uv}{L_x^1}\leq \norm{u}{L_x^2}\norm{v}{L_x^2}$, we reduce it to
showing
$$
\supl_t\norm{u}{L^2}\lesssim \norm{u}{X_{0,1/2+}}.
$$
Consider a function $f\in L^2(\rtwo\times \rone)$, such that 
$$
\widehat{u}(\xi,\tau)=f(\xi,\tau)<\tau-|\xi|^2>^{-1/2-\ve}.
$$
We have 
\begin{eqnarray*}
& &\supl_t\norm{u}{L^2_x}=
\supl_t\norm{\int\f{f(\xi,\tau)}{<\tau-|\xi|^2>^{1/2+\ve}} e^{i\tau t} d\tau}{L^2_x} \lesssim\\
& &\lesssim \norm{\norm{f(\xi,\cdot)}{L^2_{\tau}}\left(\int <\tau-|\xi|^2>^{-1-\ve}d\tau\right)^{1/2}}{L^2_{\xi}}
\lesssim \norm{f}{L^2_{\xi,\tau}}. 
\end{eqnarray*}
\end{proof}
We now turn to estimating the quintic nonlinearity $F_{quintic}$, i.e.
estimate \eqref{eq:911}.
\begin{proof}  By symmetry, it suffices to obtain  estimates only  for 
$\norm{|\nabla\beta|^2 u_1}{X_{100\ve,-1/2+2\ve}}$  in \eqref{eq:907}.
We first perform a dyadic decomposition on $\nabla\beta$ and 
$\nabla\overline{\beta}$ to get 
$$
|\nabla\beta|^2 u_1=\suml_{r_1,r_2} S_{r_1}(\nabla\beta) 
S_{r_2}(\nabla\overline{\beta}) u_1.
$$
We split into two pieces, $\max(r_1,r_2)\leq 1$ and $\max(r_1,r_2)\geq 1$.

For the small frequency case, the argument goes along the lines of the 
estimate  for $S(h,w,v)$ in \eqref{eq:2317}. Observe that
$$ S_{r_1}(\nabla\beta) 
S_{r_2}(\nabla\overline{\beta})=S_{\max(r_1,r_2)}(S_{r_1}(\nabla\beta) 
S_{r_2}(\nabla\overline{\beta})).$$ 
Therefore, just as in the estimate for $S(h,w,v)$
\begin{eqnarray*}
& &\suml_{r_1,r_2\leq 1}\norm{S_{r_1}(\nabla\beta) 
S_{r_2}(\nabla\overline{\beta})u_1}{X_{100\ve,-1/2+2\ve}}\lesssim \\
&\lesssim&\suml_{r_1,r_2\leq 1, R}\norm{S_{r_1}(\nabla\beta) 
S_{r_2}(\nabla\overline{\beta})S_R(u_1)}{X_{100\ve,-1/2+2\ve}}
\lesssim \suml_{r_1,r_2\leq 1, R} (R\max(r_1,r_2))^{2\ve}\times  \\
&\times &
\min(1,\left(\f{\max(r_1,r_2)}{R}\right)^{1/2})\norm{S_{r_1}(\nabla\beta)
S_{r_2}(\nabla\overline{\beta})}{L^2}\norm{S_R(u_1)}{X_{100\ve,1/2+\ve}}
\lesssim\\
&\lesssim &\norm{\nabla\beta}{L^4(\rtwo\times \rone)}^2
\norm{u_1}{X_{100\ve,1/2+\ve}}.
\end{eqnarray*}
To estimate $\norm{\nabla\beta}{L^4(\rtwo\times \rone)}$, one uses 
Sobolev embedding with one derivative in the space variable. We get
$$
\norm{\nabla\beta}{L^4_x(\rtwo)}\lesssim \norm{\nabla^2\beta}
{L^{4/3}_x(\rtwo)}.
$$ 
From the definition of $\beta$, the boundedness of the Riesz 
transforms and the bilinear estimate \eqref{eq:6755}, we have 
$$
\norm{\nabla^2\beta}
{L^4_tL^{4/3}_x}\lesssim 
\norm{\nabla^2\De^{-1} (u_1\overline{u_2})}{L^4_tL^{4/3}_x}
\lesssim \norm{u_1\overline{u_2}}{L^4_tL^{4/3}_x}
\lesssim \norm{u_1}{X_{\ve,1/2+\ve}}\norm{u_2}{X_{\ve,1/2+\ve}}.
$$ 
Thus 
$$
\norm{\nabla\beta}{L^4(\rtwo\times \rone)}\lesssim 
\norm{u_1}{X_{\ve,1/2+\ve}}\norm{u_2}{X_{\ve,1/2+\ve}}
$$
and we have shown \eqref{eq:907} for the quintic nonlinearity in 
 the small frequency case.

For the large frequency case,  we will have to show just as in the small
frequency case
$$
\suml_{r_1,r_2:\max(r_1,r_2)\geq 1, R}\norm{S_{r_1}(\nabla\beta) 
S_{r_2}(\nabla\overline{\beta})S_R(u_1)}{X_{100\ve,-1/2+2\ve}}
\lesssim \norm{u_1}{X_{100\ve,1/2+\ve}}^3\norm{u_2}{X_{100 \ve,1/2+\ve}}^2
$$
To verify that, following the argument  
in Proposition \ref{prop:lpk} and  \eqref{eq:2317} with 
$v=S_{r_1}(\nabla\beta) 
S_{r_2}(\nabla\overline{\beta})$, 
we will have to demonstrate some decay in $\max(r_1,r_2)$ for 
$\norm{S_{r_1}(\nabla\beta) 
S_{r_2}(\nabla\overline{\beta})}{L^2(\rtwo\times\rone)}$. 
More precisely,  we need  to show 
\begin{equation}
\label{eq:8uy}
\norm{S_{r_1}(\nabla\beta) 
S_{r_2}(\nabla\overline{\beta})}{L^2(\rtwo\times\rone)}\lesssim 
\max(r_1,r_2)^{-\si}
\norm{u_1}{X_{100\ve,1/2+\ve}}^2\norm{u_2}{X_{100\ve,1/2+\ve}}^2,
\end{equation}
for some $\si>5\ve$. 
By Cauchy-Schwartz \eqref{eq:8uy} reduces to proving
\begin{equation}
\label{eq:7681}
\norm{S_M(\nabla\beta)}{L^4(\rtwo\times\rone)}\lesssim M^{-\si}
\norm{u_1}{X_{100\ve,1/2+\ve}}\norm{u_2}{X_{100\ve,1/2+\ve}}.
\end{equation}
By Sobolev embedding performed in the spatial variable only, the boundedness 
of the Riesz transforms and the definition of $\beta$,  we get
$$
\norm{S_M(\nabla\beta)}{L^4_x(\rtwo)}\lesssim \norm{\nabla^2S_M(\beta)}
{L^{4/3}_x(\rtwo)}\lesssim \norm{S_M (u_1\overline{u_2})}{L^{4/3}_x}.
$$
Thus,
$$
\norm{S_M(\nabla\beta)}{L^4_{x,t}}\lesssim \norm{\nabla^2S_M(\beta)}
{L^{4}_tL^{4/3}_x(\rtwo)}\lesssim 
\norm{S_M (u_1\overline{u_2})}{L^4_tL^{4/3}_x},
$$
and by the bilinear estimate \eqref{eq:6755}, we get
\begin{eqnarray*}
& &\norm{S_M(\nabla\beta)}{L^4_{x,t}}\lesssim 
\max(\norm{S_M (u_1)}{X_{\ve,1/2+\ve}}\norm{u_2}{X_{\ve,1/2+\ve}},
\norm{S_M (u_2)}{X_{\ve,1/2+\ve}}\norm{u_1}{X_{\ve,1/2+\ve}})\lesssim  \\
& &\lesssim M^{-99\ve} \norm{u_1}{X_{100\ve,1/2+\ve}}\norm{u_2}
{X_{100\ve,1/2+\ve}}.
\end{eqnarray*}
\end{proof}
\vspace{.5cm}

\section{Null form: Estimates away from the diagonal }
\label{sec:10}

\vspace{.5cm}

This is the case when the ``null'' form is under control in the $L^\infty$
norm. By symmetry, it suffices to consider the case $|z|\leq |\eta|/2$.
Thus, we have the estimate
$$
\f{\dpr{\eta}{z^{\bot}}}{|\eta|^2}\lesssim \f{|z|}
{|\eta|}.
$$
Also,
\begin{equation}
\label{eq:vaj}
\La_{off diag}=\suml_R\int S_R (u_1\overline{u_2}) S_R(\tilde{G})=\suml_{M,N,R}
\int S_R((u_1)_M (\overline{u_2})_N) \tilde{G}_R d\eta d\mu,
\end{equation}
where
$$
\widehat{\tilde{G}}(\eta,\mu)=\intl_{|z|\leq |\eta|/2}
\f{\dpr{\eta}{z^{\bot}}}{|\eta|^2}\f{f(\eta-z,\mu-s)
<\eta-z>^{100\ve}}
{<\mu-s-|\eta-z|^2>^{1/2-2\ve}}
 \f{g(z,s)}{<z>^{100\ve}
<s-|z|^2>^{1/2+\ve}}dz ds.
$$
It is clear now that every term in the dyadic formula \eqref{eq:vaj} 
 can be estimated by corresponding term for the form $H$ from
Propostition \ref{prop:lpk} times $r/R$ and with $\ka=100\ve$, which makes 
the double summation in $r,R$ even easier.
We get 
$$
|\La_{offdiag}|\lesssim 
\norm{h_1}{2}\norm{h_2}{2}\norm{f}{2}\norm{g}{2}.
$$

\section{Null Form: Diagonal estimates}
\label{sec:11}

\vspace{.5cm}

In this section, we decompose the regions of the integration in \eqref{eq:100}
in such a way as to accomodate the behavior of the ``null'' form.
Let us first represent the integration region over $\eta$ and $z$
as a union  of dyadic annuli of the form  $|\eta| \sim |z|\sim R$. 
Denote $\theta=\angle(\eta,z)$. Observe that if  $|\theta|>1/100$, 
we can control
the $L^\infty$ norm of the ``null'' form in a similar manner as 
in Section \ref{sec:1}, and
that will do in that case. By symmetry, we further assume  that 
$0<\theta<1/100$.
We decompose in the angular variable $\theta$ in a dyadic manner 
as $\theta\to 0$.
Observe that
\begin{eqnarray*}
& &|\eta-z|^2=(|\eta|-|z|)^2+2|\eta||z|(1-\cos(\theta))\sim (|\eta|-|z|)^2+R^2\theta^2,\\
& &\dpr{\eta}{z^{\bot}}=|\eta||z|\sin(\theta)\sim R^2\theta.
\end{eqnarray*}
Obviously the size of $|\eta|-|z|$ is important at this stage, so we make the following  partition
of the area of integration
\begin{eqnarray*}
& &\ca_l=\left\{ (\eta,z): ||\eta|-|z||\sim 2^l R\theta   \right\},\q   l\geq 1,\\
& &\ca_0=\left\{ (\eta,z): ||\eta|-|z||\leq  R\theta   \right\}.
\end{eqnarray*}
We will concentrate on the  set $\ca_0$ and in the end we will explain
how to obtain similar estimates
when integrating on $\ca_l$ with the corresponding exponential decay in $l$.
\subsection{The ``really diagonal'' case}
To summarize, we aim at controlling the expression
\begin{eqnarray}
\label{eq:200}
& &\suml_{R,\theta}\intl_{|\eta|,|z| \sim R, ||\eta|-|z||\leq R\theta}
\f{h_1(\xi-\eta+z,\tau-\mu)<\xi-\eta+z>^{-100\ve}}{<\tau-\mu-|\xi-\eta+z|^2>^
{1/2+\ve}} \f{h_2(-\xi,-\tau)<\xi>^{-100\ve}}{<\tau-|\xi|^2>^{1/2+\ve}} \\
& &
 \f{\dpr{\eta}{z^{\bot}}}{|\eta-z|^2}
\f{f(\eta,\mu-s)<\eta>^{100\ve}}{<\mu-s-|\eta|^2>^{1/2-2\ve}}
\f{g(-z,s)<z>^{-100\ve}}{<s-|z|^2>^{1/2+\ve}}
d\xi d\eta dz d\tau d\mu ds.
\nonumber
\end{eqnarray}
Note that
\begin{eqnarray*}
& &\left|\f{\dpr{\eta}{z^{\bot}}}{|\eta-z|^2}\right|\sim \theta^{-1},\\
& &<\eta>^{100\ve}<z>^{-100\ve}\lesssim 1,
\end{eqnarray*}
whenever $(\eta,z)\in \ca_0$.

Take a  partition the annulus $|\eta|\sim R$ into a finite intersection
familly of cubes  with sides $R\theta$
$$
\{|\eta|\sim R\}=\bigcup\limits_\nu Q_\nu(\eta_\nu, R\theta).
$$
By our assumptions on $\eta$ and $z$, it is clear that
whenever $\eta\in Q_{\nu_0}$, then $z \in Q_{\nu_0}^*$, where for a cube
$Q$ we denote by $Q^*$ the cube with the same center and four times longer
sides.

The Schur's test ( cf. Lemma 3.11, \cite{Tao} ) and the finite
intersection property of the familly $\{Q_\nu\}_\nu$ imply
$$
\norm{\chi_{\{(\eta,z)\in \ca_0\}}}{[4,R^{2+1}]}\lesssim
\norm{\suml_{\nu} \chi_{\{\eta\in Q_\nu, z\in Q_\nu^*\}}}{[4,R^{2+1}]}
\lesssim \supl_{\nu}\norm{\chi_{\{\eta,z \in Q_\nu^*\}}}{[4,R^{2+1}]}.
$$
Thus, we may further restrict the region of integration in \eqref{eq:200} to a given cube
$Q_0(\eta_0,R\theta)$. By the box localization principle (Lemma \ref{le:2}),
since we have managed to restrict two of the variables to a box of size $R\theta$, we may do so for
the other two variables as well. Thus, we are lead to 
consider multipliers of the form
$\displaystyle\chi_{\{(\eta,z)\in Q_0(\eta_0,R\theta), \xi\in B(\xi_0,R\theta)\}}$, where $|\eta_0|\sim R$ and
$\xi_0\in \rtwo$. It is worth mentioning that $\eta_0$ and $\xi_0$ may depend on $R, \theta$.
Denote
\begin{eqnarray*}
L_1=\tau-\mu-|\xi-\eta+z|^2, &  L_2=\tau-|\xi|^2 ;\\
L_3=\mu-s-|\eta|^2, & L_4=s-|z|^2.
\end{eqnarray*}
To control \eqref{eq:200}, we need to show that
\begin{eqnarray}
\label{eq:300}
& &\suml_{R,\theta}
\supl_{\begin{array}{l}\eta_0,\xi_0\in \rtwo \\|\eta_0|\sim R
\end{array}} \f{\max(<|\xi_0|>,<R\theta>)^{-100\ve}}{\theta}\times\\
&\times& \norm{\f{\chi_{\{\eta,z\in Q_0(\eta_0,R\theta),|\xi-\xi_0|\leq R\theta\} }}
{(<L_1><L_2><L_4>)^{1/2+\ve}
<L_3>^{1/2-2\ve}}}{}\lesssim 1,
\nonumber
\end{eqnarray}
where we have used the fact that $<\xi-\eta+z><\xi>\geq 
\max(<|\xi_0|>,<R\theta>)$.

Next,  we write the equivalent quatrilinear
form $\La_0$  representing the multiplier in \eqref{eq:300}.
\begin{eqnarray}
\label{eq:400}
\La_0(h_1,h_2,f,g) &=&
\intl_{
\begin{array}{l}
 z,\eta\in Q_0(\eta_0,R\theta),\\
|\xi-\xi_0|< R\theta
\end{array}}
\f{h_1(\xi-\eta+z,\tau-\mu)}{<L_1>^
{1/2+\ve}} \f{h_2(-\xi,-\tau)}{<L_2>^{1/2+\ve}}  \\
&\times& \f{f(\eta,\mu-s)}{<L_3>^{1/2-2\ve}}
\f{g(-z,s)}{<L_4>^{1/2+\ve}}
d\xi d\eta dz d\tau d\mu ds.
\nonumber
\end{eqnarray}
The following lemma allows us to dramatically reduce the number of cases.
\begin{lemma}
\label{le:claim1}
With the restrictions in the integration  in \eqref{eq:400}, either $\lfour\gtrsim R^2$ or $|\xi_0|\gtrsim R/\theta$.
\end{lemma}
\begin{proof}
Assume otherwise. Then $R^2\gg |L_3+L_4|=|\mu-|\eta|^2-|z|^2|$. Since $|\eta|^2+|z|^2\sim R^2$,
it follows that $\mu\sim R^2$.

On the other hand $|\mu-|\xi|^2+|\xi-\eta+z|^2|=|L_2-L_1|\ll R^2$.
But $||\xi-\eta+z|^2-|\xi|^2|\leq |\eta-z| (|\xi|+|\xi-\eta+z|)\lesssim
R\theta\max(R\theta,|\xi_0|)\ll R^2$.
Thus $\mu\ll R^2$, a contradiction.
\end{proof}
We will actually show that
$$
\supl_{\eta_0,\xi_0:|\eta_0|\sim R}
|\La_0(h_1,h_2,f,g)|\leq C_{R,\theta,\xi_0},
$$
for suitable $C_{R,\theta,\xi_0}$, such that
\begin{equation}
\label{eq:500}
\suml_{R,\theta} \supl_{\xi_0}\f{\max(<|\xi_0|>,<R\theta>)^{-100\ve}}{\theta} C_{R,\theta,\xi_0}\lesssim 1.
\end{equation}
\vspace{.5cm}
\begin{center}
{\large\bf Case 1. $\lfour\gtrsim R^2$.}
\end{center}
\vspace{.5cm}

A subcase that can be  easily handled is when $\lthree\gtrsim R^2\theta$.
\newline
$\bullet$ {\it\bf Case 1.1.} $\lthree\gtrsim R^2\theta$
\begin{proof}
Observe that the multiplier from \eqref{eq:400} has the form
\begin{equation}
\label{eq:550}
\chi_{\left\{\begin{array}{l}
(\tau_i,\xi_i), \xi_i\in Q_i(\xi_i^0,R\theta), \\
\tau_i\pm |\xi_i|^2\sim L_i \end{array}\right\}}
\end{equation}
where $Q_i(\xi_i^0,R\theta)$ are cubes with sidelength $R\theta$.
Since all the variables are well localized, we can use \eqref{eq:11} to estimate the $[4,R^{2+1}]$
norm of the multiplier in \eqref{eq:550}.
Indeed, let us assume for simplicity that $L_1=\lone, L_2=\ltwo$. Then
choose $A, B$ in \eqref{eq:11} so that
\begin{eqnarray*}
A &=&\{(\tau_1,\xi_1):|\xi_1-\xi_1^0|\leq R\theta, \tau_1-|\xi_1|^2\sim L_1\},\\
B &=&\{(\tau_2,\xi_2):|\xi_2-\xi_2^0|\leq R\theta, \tau_2+|\xi_2|^2\sim L_2\}.
\end{eqnarray*}
For fixed $\xi_1,\xi_2$, $\tau_1$ and $\tau_2$ span intervals of length $L_1$ and $L_2$ respectively.
Therefore, since $\xi_1, \xi_2$ are both within a ball with radius $R\theta$,
we obtain from \eqref{eq:11}
\begin{equation}
\label{eq:570}
\norm{\chi_{\left\{\begin{array}{l}
(\tau_i,\xi_i), \xi_i\in Q_i(\xi_i^0,R\theta), \\
\tau_i\pm |\xi_i|^2\sim L_i \end{array}\right\}}}{[4,R^{2+1}]}\lesssim (\lone\ltwo)^{1/2}R^2\theta^2.
\end{equation}
Based on \eqref{eq:570}, we have
$$
|\La_0(h_1,h_2,f,g)\chi_{\lthree\gtrsim R^2\theta}|\lesssim \suml_{\lfour\gtrsim R^2,\lthree\gtrsim R^2\theta,\lone,\ltwo}
\f{R^2\theta^2}{(<\lone><\ltwo>)^{\ve}\lfour^{1/2-2\ve}\lthree^{1/2+\ve}}\lesssim R^{2\ve}\theta^{3/2-\ve}, 
$$
which implies \eqref{eq:500}.
\end{proof}
$\bullet${\it\bf Case 1.2.} $\lthree\lesssim R^2\theta$. \\
To avoid the enormous amount of cases to consider, we make the following
reduction. Observe that $L_3$ and $L_4$
appear symmetrically (they are both of the type $\tau_i-|\xi_i|^2$), except in the power that they
have in the denominator. Thus, since $L_3$ appear with a lesser power, it will be enough to consider the case
$L_3\geq L_4$, that is $\lfour$ does not fall into $L_4$ itself. With this 
reduction Case 1.2 breaks into five different subcases. More precisely, 
we subdivide Case 1.2 into

\begin{itemize}
\item Case 1.2.1 $\lfour=L_4$ or $L_2$.
\item Case 1.2.2 $(L_1,L_2)=(\lone,\ltwo)$ or $(\ltwo,\lone)$; 
$|\xi_0|\nsim R$ and $L_3=\lfour$.
\item Case 1.2.3 $(L_2,L_4)=(\lone,\ltwo)$ or $(\ltwo,\lone)$; 
$|\xi_0|\nsim R$ and $L_3=\lfour$.
\item Case 1.2.4 $(L_1,L_4)=(\lone,\ltwo)$ or $(\ltwo,\lone)$; 
$|\xi_0|\nsim R$ and $L_3=\lfour$.
\item Case 1.2.5 $|\xi_0|\sim R$ and then we are considering \\
 Case 1.2.5a) $\theta>1/R$ \\ 
 Case 1.2.5.b) $\theta<1/R$ since the relative sizes of $R\theta$ and 
$\xi$ will matter.
\end{itemize}

We  first dispose of the case, when $\lfour=L_1$ or
$\lfour=L_2$.\\ \\ 
{\it\bf Case 1.2.1.} $\lfour=L_1$ or $\lfour=L_2$,
$\lthree\lesssim R^2\theta$, $\lfour\gtrsim R^2$.\\
\begin{proof}
An application of the
Cauchy-Schwartz yields $$ \La_0(h_1,h_2,f,g)\lesssim
\supl_{\norm{h}{2}=1}|\La_1(h_1,h_2,h)|\supl_{\norm{h}{2}=1}|\La_2(f,g,h)|,
$$ where 
$$
\La_1(h_1,h_2,h)= \intl_{\begin{array}{l}
|\xi-\xi_0|\lesssim  R\theta, \\
|\tilde{\eta}|\lesssim R\theta
\end{array}}
\f{h_1(\xi-\tilde{\eta},\tau-\mu)}{<\tau-\mu-|\xi-\tilde{\eta}|^2>^
{1/2+\ve}} \f{h_2(-\xi,-\tau)}{<\tau-|\xi|^2>^{1/2+\ve}} h(\tilde{\eta},\mu)d\xi d\tau d\tilde{\eta}d\mu
$$
and
$$
\La_2(f,g,h)=\intl_{\begin{array}{l}
|\eta-\eta_0|\lesssim  R\theta, \\
|z-\eta_0|\lesssim R\theta
\end{array}} \f{f(\eta,\mu-s)}{<\mu-s-|\eta|^2>^{1/2-2\ve}}
\f{g(-z,s)}{<s-|z|^2>^{1/2+\ve}} h(z-\eta,-\mu)dzd\eta d\mu ds. $$
Observe that by Lemma \ref{le:5}, we can estimate
\begin{equation}
\label{eq:595} |\La_2|\lesssim \suml_{L_3, L_4 \ll \R^2\theta}
\f{(L_4)^{1/2}(R\theta/R)^{1/2}L_3^{1/2}}{<L_3>^{1/2-2\ve}<L_4>^{1/2+\ve}}
\lesssim R^{4\ve}\theta^{1/2}.
\end{equation}
By Lemma \ref{le:11}  we need to compute
$$
|H|=||\xi-\tilde{\eta}|^2-|\xi|^2+|\tilde{\eta}|^2|\lesssim
|\tilde{\eta}|(|\xi|+|\tilde{\eta}|)\lesssim R\theta \max(R\theta,|\xi_0|)
$$
Then the estimates are
\begin{equation} \label{eq:590} |\La_1|\lesssim
\suml_{L_1,L_2,\max(L_1,L_2)\gtrsim R^2}\f{\min(L_1,L_2)^{1/2}
(R\theta/\max(R\theta,|\xi_0|))^{1/2}(R\theta\max(R\theta,|\xi_0|))^{1/2}}
{(<L_1><L_2>)^{1/2+\ve}}\lesssim \f{\theta}{R^{2\ve}}.
\end{equation}
Combining \eqref{eq:595} and \eqref{eq:590} gives \eqref{eq:500}
in Case 1.2.1.
\end{proof}
We will postpone the somewhat 
peculiar case $|\xi_0|\sim R$ for later on.

We consider the case, where $L_1$ and $L_2$ are the two smallest numbers in the sequence
$L_1, L_2, L_3, L_4$.\\ \\
{\it\bf Case 1.2.2} $(L_1,L_2)=(\lone,\ltwo)$ or 
$(\ltwo,\lone), |\xi_0|\nsim R $, $L_3=\lfour\gtrsim R^2$.
\begin{proof}
In that case, we will fully use the quatrilinear form, instead of
relying on Cauchy-Schwartz and then deal with the resulting
trilinear forms. We estimate the multiplier in \eqref{eq:400} by
\eqref{eq:11}. We have an upper bound of 
\begin{equation}
\label{eq:601}
 \suml_{L_1,L_2,L_3,L_4}\f{\left|
\left\{\{(\tau_1,\xi_1),(\tau_2,\xi_2)\}\in(\Om_1\times\Om_2):(\tilde{\tau}-\tau_1-\tau_2,
\tilde{\xi}-\xi_1-\xi_2)\in
\Om_4\right\}\right|^{1/2}}{(<L_1><L_2><L_4>)^{1/2+\ve}<L_3>^{1/2-2\ve}},
\end{equation}
where $\tilde{\tau},\tilde{\xi}$ are fixed and $$\Om_1=\left\{
(\tau_1,\xi_1):\tau_1-|\xi_1|^2\sim L_1,\qq \xi_1=\xi-\eta+z,
|\xi-\xi_0|\leq R\theta,\q |\eta-z|\leq R\theta \right\} $$
$$\Om_2=\left\{ (\tau_2,\xi_2):\tau_2+|\xi_2|^2\sim L_2,\qq
\xi_2=-\xi,\q |\xi-\xi_0|\leq R\theta\right\} $$ $$\Om_4=\left\{
(\tau_4,\xi_4):\tau_4-|\xi_4|^2\sim L_4, \qq \xi_4=-z,\q
|z-z_0|\leq R\theta\right\} $$ Note that for a fixed spatial
variables the time variables span intervals of length $L_1$ and
$L_2$ respectively. Also, we have
\begin{equation}
\label{eq:604}
\tilde{\tau}-|\tilde{\xi}-\xi_1-\xi_2|^2-|\xi_1|^2+|\xi_2|^2=L_1+L_2+L_4.
\end{equation}
For a {\bf fixed} $\xi_2$, we have (based on \eqref{eq:604})
$$
 |\tilde{\xi}-\xi_1-\xi_2|^2 +|\xi_1|^2=const
+O(\lthree).$$ By the parallelogram law,
\begin{equation}
\label{eq:602} |(\tilde{\xi}-\xi_2)/2-\xi_1|^2=const + O(\lthree)
\end{equation}

Furthermore, since $\tilde{\xi}-(\xi_1+\xi_2)\in B(z_0, R\theta)$
and $|\xi_1+\xi_2|\lesssim R\theta$, we infer $\tilde{\xi}\sim R$.
Thus taking into account that $|\xi_0|\nsim R$, we conclude that
$(\tilde{\xi}-\xi_2)/2-\xi_1\sim \max(R,|\xi_0|)$ and therefore by
\eqref{eq:602}, $\xi_1$ is contained in an annulus with radius
$\max(R,|\xi_0|)$ and thickness $\lthree/\max(R,|\xi_0|)$. Observe
that $\xi_1$ is also in a ball with radius $R\theta$,
therefore it belongs to a rectangle with sides $R\theta$ and
$\lthree/\max(R,|\xi_0|)$. Finally, since $\xi_2$ belongs to a
ball with radius $R\theta$, one estimates \eqref{eq:601} by $$
\suml_{L_1,L_2, L_3, L_4}\f{(L_1 L_2)^{1/2}R\theta
\left(\f{\lthree
R\theta}{\max(R,|\xi_0|)}\right)^{1/2}}{(<L_1><L_2><L_4>)^{1/2+\ve}L_3^{1/2-2\ve}}
\lesssim R^{4\ve}\theta^{3/2}, $$
thus implying \eqref{eq:500}.
\end{proof}
{\it\bf Case 1.2.3} $(L_2,L_4)=(\lone,\ltwo)$ or 
$(\ltwo,\lone),\q |\xi_0|\nsim R$, $L_3=\lfour$.

\begin{proof}
This case is very similar to Case 1.2.2. We estimate the
multiplier by
\begin{equation}
\label{eq:610}
 \suml_{L_1,L_2,L_3,L_4}\f{\left|
\left\{\{(\tau_1,\xi_1),(\tau_2,\xi_2)\}\in(\Om_2\times\Om_4):(\tilde{\tau}-\tau_1-\tau_2,
\tilde{\xi}-\xi_1-\xi_2)\in
\Om_1\right\}\right|^{1/2}}{(<L_1><L_2><L_4>)^{1/2+\ve}<L_3>^{1/2-2\ve}},
\end{equation}
where $\Om_1, \Om_2, \Om_4$ are the sets defined before. We have $$
\tilde{\tau}-|\tilde{\xi}-\xi_1-\xi_2|^2-|\xi_2|^2+|\xi_1|^2=L_1+L_2+L_4=O(\lthree)$$ 
For {\bf fixed} $\xi_1$ , we have by the parallelogram law
\begin{equation}
\label{eq:671}
 |(\tilde{\xi}-\xi_1)/2-\xi_2|=const+O(\lthree).
\end{equation}
Since $\tilde{\xi}-\xi_1-\xi_2\in B(\xi_0, R\theta)$ and
$|\xi_1+\xi_0|\lesssim R\theta$, it follows that
$|\tilde{\xi}-\xi_2|\lesssim R\theta$. In particular, after taking into
account that $\xi_0\nsim R$, we obtain
$|(\tilde{\xi}-\xi_1)/2-\xi_2|\sim \max(R,|\xi_0|).$ By
\eqref{eq:671}, one has that $\xi_2$ is contained in an annulus
with radius $\max(R,|\xi_0|)$ and thickness
$\lthree/\max(R,|\xi_0|)$. Since $\xi_2$ is also contained in a
ball with radius $\R \theta$, we have that $\xi_2$ is contained in
a rectangle with sidelengths $R\theta$ and
$\lthree/\max(R,|\xi_0|)$ for every fixed $\xi_1$. The usual
observation that $\tau_1$ and $\tau_2$ sweep intervals of length
$L_2$ and $L_4$ respectively, leads us to estimate \eqref{eq:610} by
$$ \suml_{L_1,L_2,L_3,L_4}\f{(<L_2><L_4>)^{1/2} R\theta \left(\f{\lthree
R\theta}{\max(R,|\xi_0|)}\right)^{1/2}}{(<\lone><\ltwo><\lthree>)^{1/2+\ve}\lfour^{1/2-2\ve}}
\lesssim R^{4\ve}\theta^{3/2},$$
which again implies \eqref{eq:500}.
\end{proof}

{\it\bf Case 1.2.4} $(L_1,L_4)=(\lone,\ltwo)$ or 
$(\ltwo,\lone), |\xi_0|\nsim R$, $L_3=\lfour\gtrsim R^2$.

\begin{proof}
We estimate the norm of the multiplier by
\begin{equation}
\label{eq:701}
 \suml_{L_1,L_2,L_3,L_4}\f{\left|
\left\{\{(\tau_1,\xi_1),(\tau_2,\xi_2)\}\in(\Om_1\times\Om_4):(\tilde{\tau}-\tau_1-\tau_2,
\tilde{\xi}-\xi_1-\xi_2)\in
\Om_2\right\}\right|^{1/2}}{(<L_1><L_2><L_4>)^{1/2+\ve}<L_3>^{1/2-2\ve}},
\end{equation}
where $\Om_1, \Om_2, \Om_4$ are the sets defined in Case 1.2.2.
Like in the previous cases, we have  a relation involving some of
the variables. Here, we have $$
\tilde{\tau}-|\xi_1|^2-|\xi_2|^2+|\tilde{\xi}-\xi_1-\xi_2|^2=L_1+L_2+L_4=O(\lthree).
$$ We change variables $\la_1=\xi_1+\xi_2$, $\la_2=\xi_1-\xi_2$
and we are interested in the measure of the corresponding set in
\eqref{eq:701}. Since the Jacobian of the transformation  is two,
we pass to the new variables. Fix $\la_1$. Observe also that since
$|\xi_0|\nsim R$, $|\la_2|\sim \max(R,|\xi_0|)$.  We have then by
the parallelogram law
\begin{eqnarray*}
|\xi_1|^2+|\xi_2|^2 &=& const+O(\lthree), \\ |\xi_1-\xi_2|^2 & =&
const+O(\lthree).
\end{eqnarray*}
That implies that for fixed $\la_1$, $\la_2$ is contained in an
annulus with thickness $\lthree/\max(R,|\xi_0|)$. On the other
hand, $\la_2$ is contained in a ball with radius $R\theta$. These
estimates, together with the usual observations that $\tau_1,
\tau_2$ are in a intervals of length $L_1$ and $L_4$ respectively,
and the fact that $\la_1$ sweeps a ball with radius $R\theta$
imply the following bound for \eqref{eq:701} $$ \suml_{L_1, L_2,
L_3, L_4} \f{(L_1L_4)^{1/2} R\theta \left(\f{\lthree
R\theta}{\max(R,|\xi_0|)}\right)^{1/2}}{(<\lone><\ltwo><\lthree>)^{1/2+\ve}\lfour^{1/2-2\ve}}
\lesssim R^{4\ve}\theta^{3/2}, $$
which clearly implies \eqref{eq:500}.
\end{proof}
Finally, we deal with the  case  $|\xi_0|\sim R$.\\\\
$\bullet${\it\bf Case 1.2.5} 
$|\xi_0|\sim R, \lfour\gtrsim R^2, \lthree\lesssim  R^2\theta$ \newline
In that case, the relative size of $R\theta$ and $\xi_0$ will mater, so we will split into two subcases.\\\\
{\it\bf Case 1.2.5a)} $\theta>1/R$\\
{\bf Remark} This case is vacuous if $R<1$.
\begin{proof}
Apply Cauchy-Schwartz to $\La_0$ obtain
$$
|\La_0(h_1,h_2,f,g)|\lesssim \supl_{\norm{h}{2}=1}|\La_1(h_1,h_2,h)|\supl_{\norm{h}{2}=1}|\La_2(f,g,h)|,
$$
where
$\La_1$ and $\La_2$ were defined in Case 1.2.
Compute again  $H=||\tilde{\eta}|^2+|\xi-\tilde{\eta}|^2-|\xi|^2|\lesssim R^2\theta$.
Thus, based on  the estimates  in Lemma \ref{le:11}, we conclude that
 $$ |\La_1|\lesssim
\suml_{L_1,L_2}\min(L_1, L_2)^{1/2}(R\theta/R)^{1/2}\f{\min(R^2\theta,\max(L_1,L_2)/\theta)^{1/2}}{(<L_1><L_2>)^{1/2+\ve}}
$$ and by the estimate for $m_1$ in Lemma \ref{le:5} $$
|\La_2|\lesssim\suml_{L_3,L_4}
\f{\min(L_3,L_4)^{1/2}(R\theta/R)^{1/2}\min(\max(L_3,L_4),R^2\theta)^{1/2}}{<L_3>^{1/2-2\ve}<L_4>^{1/2+\ve}}.
$$ 
We have the estimate (after quickly going through the appropriate cases - $\lfour$ is either $L_1$ or $L_2$ or $\lfour$ 
is either $L_3$ or $L_4$)
$$|\La_1||\La_2|\lesssim
R^{4\ve}\theta.$$ 
Thus we can sum up the expression in \eqref{eq:500} as follows
$$
\suml_{R}\q\suml_{\theta>1/R, \theta\textup{dyadic}} \f{R^{-100\ve}}{\theta} R^{4\ve}\theta\lesssim \suml_{R} 
R^{-96\ve}\ln{R}\lesssim 1
$$
\end{proof}
{\it\bf Case 1.2.5b)} $\theta<1/R$. 
\begin{proof}
We concentrate on the high frequency case $R>1$. The case $R<1$ is 
trivial, because then $<\lone>=<\ltwo>=<\lthree>=1$ and one easily 
estimates (see estimates below).
For simplicity, we assume once again that $L_3=\lfour\gtrsim R^2$. Observe that $L_3$ appears with the smallest
power in the denominator and that should be the worst case for the maximum to occur. Moreover later in 
the proof regarding that case, we will see that we could perform the same argument with any other configuration
of $\lone, \ltwo, \lthree, \lfour$.
We use again the quatrilinear form $\La_0$. By \eqref{eq:11}, we  get
\begin{equation}
\label{eq:567}
\suml_{L_1,L_2,L_3,L_4}\f{\left|
\left\{\{(\tau_1,\xi_1),(\tau_2,\xi_2)\}\in(\Om_1\times\Om_4):(\tilde{\tau}-\tau_1-\tau_2,
\tilde{\xi}-\xi_1-\xi_2)\in
\Om_2\right\}\right|^{1/2}}{(<L_1><L_2><L_4>)^{1/2+\ve}<L_3>^{1/2-2\ve}}.
\end{equation}
We have the relation 
\begin{equation}
\label{eq:950}
\tilde{\tau}+|\tilde{\xi}-\xi_1-\xi_2|^2-|\xi_1|^2-|\xi_2|^2=L_1+L_2+L_4=O(\lthree).
\end{equation}
There are two distinct possibilities now. Either $|\xi_0+z_0|\geq R$ or $|\xi_0-z_0|\geq R$ (or both). We
show  the desired estimate, for the case $|\xi_0+z_0|\geq R$, the other case being similar.
Observe that $|\xi_1-\xi_2|=|\xi_0+z_0|+O(R\theta)\gtrsim R$. We introduce again the 
new variables $\la_1=\xi_1+\xi_2, \la_2=\xi_1-\xi_2$ and we fix $\la_1$. The parallelogram law and \eqref{eq:950} imply
\begin{eqnarray*}
|\xi_1|^2+|\xi_2|^2 &=& const+O(\lthree),\\
|\xi_1-\xi_2|^2 &=& const+O(\lthree)
\end{eqnarray*}
and thus, we have that for fixed $\la_1$, $\la_2$ is contained in an annulus with thickness $\lthree/R$. On the 
other hand it is contained in a ball with radius $R\theta$. The usual observation that $\tau_1, \tau_2$ span 
intervals of length $L_1, L_4$, gives us the following estimate for \eqref{eq:567}
\begin{eqnarray*}
& &\suml_{L_3\gtrsim R^2, \lthree\ll R^2\theta}\f{(L_1L_4)^{1/2}R\theta( R\theta \lthree/R)^{1/2}}
{(<\lone\ltwo\lthree>)^{1/2+\ve} <L_3>^{1/2-2\ve}}\lesssim 
\suml_{\lthree\ll R^2\theta}\theta (\lthree\theta)^{1/2} R^{4\ve}\lesssim R\theta^2 R^{4\ve},
\end{eqnarray*}
which implies 
$$
\suml_{R} \suml_{\theta<1/R} \f{R^{-100\ve} R\theta^2 R^{4\ve}}{\theta}\lesssim 1,
$$
which is the desired estimate \eqref{eq:500}.
\end{proof}

Now, we pass to the other possibility alluded to in Lemma \ref{le:claim1}, 
namely
that $|\xi_0|\gtrsim R/\theta$. Since we have exhausted the cases, where 
$\lfour \gtrsim R^2$, we will consider only $\lfour \ll R^2$.
\vspace{.5cm}
\begin{center}
{\large\bf Case 2.} $|\xi_0|>R/\theta$, $\lfour\ll R^2$.
\end{center}
\vspace{0.5cm}
\begin{proof}
We will have to apply Cauchy-Schwartz's inequality with a reorganized pairs of functions. We do that in order to  take advantage of the 
disparity in the sizes of $|\xi_0|$ and $|\eta_0|$.
$$
\La_0(h_1,h_2,f,g)\lesssim
\supl_{\norm{h}{2}=1}|\La_1(h_1,h_2,h)|\supl_{\norm{h}{2}=1}|\La_2(f,g,h)|,
$$ where 
\begin{eqnarray*}
\La_1(h_1,h_2,h) &=& \intl_{\begin{array}{l}
|\xi-\xi_0|\lesssim  R\theta, \\
|\eta-\eta_0|\lesssim R\theta \\
\lfour\ll R^2
\end{array}}
\f{h_1(\xi,\tau_1)}{<\tau_1-|\xi|^2>^
{1/2+\ve}} \f{f(\eta,\tau_2)}{<\tau_2-|\eta|^2>^{1/2-2\ve}} \times \\
&\times&h(-\xi-\eta,-\tau_1-\tau_2)d\xi  d\eta d\tau_1 d\tau_2
\end{eqnarray*}
and
\begin{eqnarray*}
\La_2(f,g,h) &=& \intl_{\begin{array}{l}
|z-\eta_0|\lesssim  R\theta, \\
|\xi-\xi_0|\lesssim R\theta
\end{array}} \f{h_2(-\xi,-\tau_1)}{<\tau_1-|\xi|^2>^{1/2+\ve}}
\f{g(z,\tau_2)}{<\tau_2-|z|^2>^{1/2+\ve}} \times \\
&\times&h(\xi-z,-\tau_1-\tau_2)
d\xi dz d\tau_1 d\tau_2. 
\end{eqnarray*}
For $\La_1$, we are in a position to use Lemma \ref{le:5}. Since 
$\lfour\ll R^2$, we  have
$$
|\La_1|\lesssim \suml_{L_1,L_3\ll R^2}
\f{\min(L_1,L_3)^{1/2}}{L_1^{1/2+\ve}L_3^{1/2-2\ve}}
\left(\f{R}{R/\theta}\right)^{1/2}
\max(L_1,L_3)^{1/2}\lesssim 
R^{4\ve}\theta^{1/2}.
$$
For $\La_2$, we  use Lemma \ref{le:3} to infer
$$
|\La_2|\lesssim \theta^{1/2}.
$$
For $R>1$, we combine the estimates for 
$\La_1$ and $\La_2$ to show 
$$
\suml_{R,\theta} \left(\f{\theta}{R}\right)^{100\ve} 
\f{R^{4\ve}\theta^{1/2}\theta^{1/2}}{\theta}\lesssim 1,
$$
which is the desired inequality  \eqref{eq:500}. 
For $R<1$, the estimates above can be improved greatly and thus one estimates
in that case as well.
\end{proof}

\vspace{.5cm}

\subsection{Null form: The not so diagonal case}

\vspace{.5cm}

This is the case where the integration in the definition of $\La_0$ is over the set $\ca_l$. Note first,
that if $(\eta,z)\in \ca_l$, then 
\begin{equation}
\label{eq:987}
\left|\f{\dpr{\eta}{z^{\bot}}}{|\eta-z|^2}\right|\sim \f{1}{2^{2l}\theta}.
\end{equation}
Thus, if $2^l\theta\gtrsim 2^{-l/2}$, the $L^\infty$ norm of the ``null'' form is under control 
(with exponential decay in 
$l$) and we can  estimate as in the off-diagonal case. 
So assume from now on that $2^l\theta\ll 2^{-l/2}$. Denote 
\begin{eqnarray*}
\La_l(h_1,h_2,f,g) &=&
\intl_{(\eta,z)\in \ca_l, |\eta|,|z| \sim R}
\f{h_1(\xi-\eta+z,\tau-\mu)}{<L_1>^
{1/2+\ve}} \f{h_2(-\xi,-\tau)}{<L_2>^{1/2+\ve}}  \\
&\times& \f{f(\eta,\mu-s)}{<L_3>^{1/2-2\ve}}
\f{g(-z,s)}{<L_4>^{1/2+\ve}}
d\xi d\eta dz d\tau d\mu ds.
\end{eqnarray*}
Taking into account \eqref{eq:987},  we need to show that 
\begin{equation}
\label{eq:988}
\suml_{R, \theta, l:2^l\theta\ll2^{-l/2}} \supl_{\xi_0}\f{\max(|\xi_0|, <2^lR\theta>)^{-100\ve}}{2^{2l}\theta}|\La_l|\lesssim 1. 
\end{equation}
Partition  the annulus 
$$
\{|\eta|\sim R\}=\bigcup\limits_{\nu\in\Th} Q(\eta_\nu, R\theta),
$$
into a finite intersection familly of cubes of sidelength $R\theta$. 
Partition the set $\ca_l$ accordingly
$$
\ca_l=\bigcup\limits_{j=-2^{l-1}}^{2^{l+1}} \{(\eta,z)\in \ca_l: ||\eta|-|z|-(2^l+j)R\theta|\leq R\theta\}=
\bigcup\limits_{j} \ca_l^j.
$$
For every fixed $j$, there is a selector map $m_j:\Th\to \Th$, so that if $(\eta,z)\in \ca_l^j$ and
whenever $\eta\in Q(\eta_\nu, R\theta)$, then $z\in Q(\eta_{m_j(\nu)}, R\theta)$ and $\eta_\nu\| \eta_{m_j(\nu)}$. 
This is possible since $\theta=\angle(\eta,z)\ll 1$. Thus, by the Schur's
 test
\begin{equation}
\label{eq:999}
\norm{\La_l}{[4,R^{2+1}]}\lesssim \suml_{j=-2^{l-1}}^{2^{l+1}}\norm{\La_l^j}{[4,R^{2+1}]}
\lesssim 2^l\supl_{j} \norm{\La_l^j}{[4,R^{2+1}]},
\end{equation}
where $\La_l^j$ are in the form $\La_0$ with the integration taken in the corresponding region $\ca_l^j$. 
By the localization principle, applied to each one of the forms 
$\La_l^j$, we can further restrict the $\xi$ integration
to a ball with center $\xi_0$ and radius $R\theta$. Thus, we are lead to estimate quatrilinear forms 
\begin{eqnarray}
\label{eq:k1}
\La(h_1,h_2,f,g) &=&
\intl_{\begin{array}{l}
|\eta-\eta_0|\leq R\theta, |z-z_0|\leq R\theta, \\
|\xi-\xi_0|\leq R\theta, |\eta_0-z_0|\sim 2^l R\theta \\
\eta_0\|z_0; |\eta_0|, |z_0|\sim R
\end{array}}
\f{h_1(\xi-\eta+z,\tau-\mu)}{<L_1>^
{1/2+\ve}} \f{h_2(-\xi,-\tau)}{<L_2>^{1/2+\ve}}  \\
\nonumber
&\times& \f{f(\eta,\mu-s)}{<L_3>^{1/2-2\ve}}
\f{g(-z,s)}{<L_4>^{1/2+\ve}}
d\xi d\eta dz d\tau d\mu ds.
\end{eqnarray}
Thus to show \eqref{eq:988}, it will suffice to obtain an estimate
$$
|\La(h_1, h_2, f,g)|\leq C_{l,R, \theta}\norm{h_1}{2}\norm{h_2}{2}\norm{f}{2}\norm{g}{2},
$$
where $C_{l, R, \theta}$ satisfies
\begin{equation}
\label{eq:2000}
\suml_{R, \theta, l:2^l\theta\ll 2^{-l/2}} \supl_{\xi_0}\f{\max(|\xi_0|, <2^lR\theta>)^{-100\ve}}{2^{l}\theta}C_{l, R, \theta}\lesssim 1. 
\end{equation}
We start reviewing the proof that we gave for the boundedness of the similar
quatrilinear form $\La_0$ given in \eqref{eq:400}. First, observe 
that Lemma \ref{le:claim1} has to read now
\begin{lemma}
\label{le:claim2}
$\lfour\gtrsim R^2$ or $|\xi_0|\gtrsim R/(2^l\theta)$.
\end{lemma}
In the Case 1.1, we will have the exact same estimate regardless of
the new restriction $|\eta_0-\xi_0|\sim 2^l R\theta$, which will enable us
to add up in \eqref{eq:2000} thanks to the exponential factor in the
denominator.

In Case 1.2.1, we obtain the estimate $|\La|\lesssim R^{2\ve} (2^l\theta)^{3/2}$, rather than
$|\La|\lesssim R^{2\ve} (\theta)^{3/2}$, but that still implies the validity of \eqref{eq:2000}, since
$$
\suml_{R,\theta, l: 2^l \theta\ll 2^{-l/2}} 
\f{(2^l\theta)^{3/2}R^{2\ve}}{2^l\theta (2^lR\theta)^{100\ve}}\lesssim 1.
$$

In Cases 1.2.2, 1.2.3, 1.2.4 we will have absolutely no change in 
the estimates, hence we can add up (in $l$) in \eqref{eq:2000}, 
due to the exponential factor in the denominator.

The restrictions in Case 1.2.5a) should be changed now to 
$|\xi_0|\sim R, \lfour\gtrsim R^2, \lthree\lesssim 
 R^2\theta, 2^l\theta>1/R$. The case $R\leq 1$ is easy to estimate. Indeed,
one can proceed as in Case 1.1 to estimate the multiplier norm in 
\eqref{eq:k1} by $R^2\theta^2$ and therefore \eqref{eq:2000} follows by
$$
\suml_{l,R\leq 1, \theta\leq 1}\f{R^2\theta^2}{2^l\theta}\lesssim 1.
$$
For $R\>1$, a careful inspection of the argument, 
shows that one has an estimate 
$|\La|\lesssim 2^l\theta R^{4\ve}$, but we can still add up in 
\eqref{eq:2000}, since
$$
\suml_{R,l,\theta: 2^l\theta\ll2^{-l/2}} \f{R^{-100\ve}R^{4\ve}2^l\theta}
{2^l\theta}\lesssim \suml_{R} R^{-96\ve} 
\suml_{l,\theta: 1\leq l \lesssim \ln(R),\theta\gtrsim R^{-3}}\lesssim 
\suml_R \ln^2(R) R^{-96\ve}\lesssim 1.
$$

In Case 1.2.5 b), we treat $|\xi_0|\sim R, \lfour\gtrsim R^2, 
\lthree\lesssim R^2\theta, 2^l\theta<1/R$. The case $R<1$ can be performed as 
in the case 1.2.5a) above. For $R>1$, the argument in Case 1.2.5b) 
 shows that the
same estimate holds ( regardless of the new restriction on $\eta_0,\xi_0$). 
Thus we add up with the help of the exponential factor in \eqref{eq:2000}.

In Case 2 according to Lemma \ref{le:claim2}, 
we deal with $|\xi_0|\gtrsim R/(2^l\theta)$, 
instead of $|\xi_0|\gtrsim R/\theta$. 
For Case 2.1, the estimate is $|\La|\lesssim R^{4\ve}2^{l/2}\theta^{1+4\ve}$. Thus \eqref{eq:2000}
is bounded by
$$
\suml_{R,\theta,l: 2^l\theta\ll 1} \f{R^{4\ve}\theta 2^{l/2}}{2^l\theta <R/(2^l\theta)>^{100\ve}}
\lesssim \suml_{\theta,l:2^l\theta\ll 1} \f{(2^l\theta)^{100\ve}}{2^{l/2}}\lesssim 1.
$$

Finally, in Case 2, the restrictions are  $|\xi_0|\gtrsim R/(2^l\theta)$, 
$\lfour\ll R^2\theta^2$.
We obtain an estimate $|\La|\lesssim (2^{l}\theta)^{3/2}R^{6\ve}$, 
which is 
handled as in the 
Case 1.2.1.

\vspace{.5cm}

\section{Regularity results}
\label{sec:reg}

\vspace{.5cm}

In this section, we will show that once we have a local existence result
for data $u_0\in H^{100\ve}$ with lifespan for the solution 
$T=T(\norm{u_0}{H^{100\ve}})$, then we have local existence for 
data $u_0\in H^{k}$ with a lifespan for the solution 
$T_k=T(\norm{u_0}{H^{100\ve}},k)$.  More precisely, 
we have

\begin{theorem}(Regularity estimates)

For a given data $u_0\in H^k$, the system 
\eqref{eq:1} has an unique solution $u$ defined at least for time  
$T=T(\norm{u_0}{H^{100\ve}},k)$ and there exists a constant $C_{\ve,k}$, 
so that  
$$
\norm{u}{X_{k,1/2+\ve}}\leq C_{\ve,k}\norm{u_0}{H^k}.
$$
\end{theorem}
\begin{proof}
For our purposes it will suffice to check the statement of the theorem for
some specific sequence $(k_n)$, so that $k_n\to\infty$, since for every 
$u_0\in H^k$, we will find $n$ so that $k_{n+1}k\geq k_n$ and the solution 
has a lifespan at least $T_{k_n}$. One could obtain estimates for the
indeces in between by the Leibnitz rule for fractional differentiation
in the $X_{s,b}$ spaces. We do not pursue these however since they are not
necessary for our purposes.  

Due to the nature 
of the estimates, it will
be convenient to take $k=n+100\ve$.
We will show the theorem for a cubic nonlinearity, since the others 
are treated in the same way. The common between them is the (anti) linearity
structure that all of them exhibit. 
Take $k=1+100\ve$. Since $u_0\in H^1\subset H^{100\ve}$, it is clear that a solution 
exists and it satisfies the integral equation \eqref{eq:3}. Differentiating
 \eqref{eq:3} yields
$$
\nabla_x u(x,t)=e^{i t \De}u_0+\intl_0^t e^{i (t-\tau)\De}
\nabla_x F_u(\tau, \cdot)d\tau.
$$
Assume now that $F_u=F_{cubic}$ for simplicity.  Note that by the linearity
of $F_{cubic}$ and  the product rule
$$
\nabla_x F_{cubic}(u_1,u_2,u_3)=
F_{cubic}(\nabla u_1,u_2,u_3)+F_{cubic}( u_1,\nabla u_2,u_3)+
F_{cubic}( u_1, u_2,\nabla u_3),
$$
where $u_3=u_1$ or $u_2$ as usual.
Thus, by combining estimates \eqref{eq:2} and \eqref{eq:907}, we obtain
\begin{eqnarray*}
& &\norm{u}{X_{1+100\ve,1/2+\ve}}\lesssim 
\norm{u_0}{H^{1+100\ve}}+ 
T^\ve\norm{u_1}{X_{100\ve,1/2+\ve}}\norm{u_2}{X_{100\ve,1/2+\ve}}
\norm{\nabla u_3}{X_{100\ve,1/2+\ve}}+\ldots \\
&\lesssim& \norm{u_0}{H^{1+100\ve}}+ T^{\ve}\norm{u}{X_{1+100\ve,1/2+\ve}}
\norm{u}{X_{100\ve,1/2+\ve}}^2.
\end{eqnarray*}
Thus, taking 
$$
T^\ve \norm{u}{X_{100\ve,1/2+\ve}}^2\lesssim 
T^\ve\norm{u_0}{H^{100\ve}}^2\ll 1
$$ 
will allow us 
to hide $T^{\ve}\norm{u}{X_{1+100\ve,1/2+\ve}}
\norm{u}{X_{100\ve,1/2+\ve}}^2$ and we get (for a time 
$T=T(\norm{u}{H^{100\ve}})$)
$$
\norm{u}{X_{1+100\ve,1/2+\ve}}\lesssim 
\norm{u_0}{H^{1+100\ve}}.
$$
We can basically iterate this result in the following manner.
For $k=2+100\ve$, we proceed as follows.
Differentiate \eqref{eq:3} once more. Then, we have two groups of terms. 
In the first, the derivatives fall on different $u$, for example
$F_{cubic}(\nabla u_1, \nabla u_2, u_3)$, for the second we will have terms 
like $F_{cubic}(\nabla^2 u_1, u_2, u_3)$. All the terms  
are estimated by \eqref{eq:907}. We get
\begin{eqnarray}
\label{eq:212}
& &\norm{u}{X_{2+100\ve,1/2+\ve}}\lesssim 
\norm{u_0}{H^{2+100\ve}}+ 
T^\ve\norm{u_1}{X_{100\ve,1/2+\ve}}\norm{u_2}{X_{100\ve,1/2+\ve}}
\norm{\nabla^2 u_3}{X_{100\ve,1/2+\ve}}+\ldots \\
\nonumber
&+& T^\ve\norm{\nabla u_1}{X_{100\ve,1/2+\ve}}
\norm{\nabla u_2}{X_{100\ve,1/2+\ve}}
\norm{ u_3}{X_{100\ve,1/2+\ve}}+\ldots \\
\nonumber
&\lesssim& \norm{u_0}{H^{2+100\ve}}+ T^{\ve}\norm{u}{X_{2+100\ve,1/2+\ve}}
\norm{u}{X_{100\ve,1/2+\ve}}^2+ \ldots\\
\nonumber
&+& T^\ve\norm{\nabla u_1}{X_{100\ve,1/2+\ve}}
\norm{\nabla u_2}{X_{100\ve,1/2+\ve}}
\norm{ u_3}{X_{100\ve,1/2+\ve}}+\ldots
\end{eqnarray}
We further have
\begin{eqnarray*}
& &\norm{\nabla u_1}{X_{100\ve,1/2+\ve}}
\norm{\nabla u_2}{X_{100\ve,1/2+\ve}}
\lesssim \suml_{M,N}
 \norm{S_M \nabla u_1}{X_{100\ve,1/2+\ve}}
\norm{S_N \nabla u_2}{X_{100\ve,1/2+\ve}}
 \lesssim \\
&\lesssim& \suml_{M\geq N} 
\f{N}{M}\norm{S_M u}{X_{2+100\ve,1/2+\ve}}
\norm{S_N u}{X_{100\ve,1/2+\ve}}\lesssim 
\norm{u}{X_{2+100\ve,1/2+\ve}}\norm{u}{X_{100\ve,1/2+\ve}}.
\end{eqnarray*}
Inserting the estimates in \eqref{eq:212}, we get
$$
\norm{u}{X_{2+100\ve,1/2+\ve}}\lesssim 
\norm{u_0}{H^{2+100\ve}}+ T^{\ve}\norm{u}{X_{2+100\ve,1/2+\ve}}
\norm{u}{X_{100\ve,1/2+\ve}}^2
$$
Choosing $T$ with $T^{\ve}\norm{u}{X_{100\ve,1/2+\ve}}^2\lesssim 
T^{\ve}\norm{u_0}{H^{100\ve}}^2\ll 1$ allows 
us to hide again. We obtain
$$
\norm{u}{X_{2+100\ve,1/2+\ve}}\lesssim \norm{u_0}{H^{2+100\ve}}.
$$
In a similar fashion, one obtains the estimates for every $k=100\ve+n$. 
We omit the details.
\end{proof}

\end{document}